\documentclass[11pt,a4paper]{amsart}

\usepackage{amscd} 

\newtheorem{theorem}{Theorem}[section]
\newtheorem{lemma}[theorem]{Lemma}
\newtheorem{prop}[theorem]{Proposition}
\newtheorem{cor}[theorem]{Corollary}
\theoremstyle{definition}
\newtheorem{definition}[theorem]{Definition}
\theoremstyle{remark}
\newtheorem{remark}[theorem]{Remark}

\def\C{{\mathbb{C}}} \def\Z{{\mathbb{Z}}} \def\N{{\mathbb{N}}} 
\def\P{{\mathbb{P}}} \def\A{{\mathbb{A}}} \def\F{{\mathbb{F}}} 
\def\G{{\mathbb{G}}} \def\O{{\mathcal O}} \def\Q{{\mathbb{Q}}}

\def\Cal{\mathcal} 

\def\Sl{{\rm Sl}} \def\Gl{{\rm Gl}} \def\Aut{{\rm Aut}}
\newcommand{\Gal}{\operatorname{Gal}}
\newcommand{\End}{\operatorname{End}}
\newcommand{\Hom}{\operatorname{Hom}}

\def\d{\partial} \def\s{\sigma} \def\t{\tau}
\def\ph{\varphi} \def\g{\gamma}

\def\wh{\widehat} \def\bs{\backslash}
\def\Spec{{\rm Spec}}
\def\lr{\longrightarrow} \def\ll{\longleftarrow}
\def\td{\ph^{\rm td}} \def\tdm{\ph^{\rm td}_{\mf m}}
\def\ld{\lambda^{\rm td}} \def\ldm{\lambda^{\rm td}_{\mf m}}

\newcommand{\ov}[1]{\overline{#1}}
\newcommand{\mf}[1]{\mathfrak{#1}}

\title{Drinfeld modular curve and Weil pairing}

\author{Gert-Jan van der Heiden}
\address{
Department of Mathematics R{\sl u}G \\
P.O. Box 800 \\
9700 AV Groningen \\
the Netherlands
}
\email{gertjan@math.rug.nl}

\date\today

\begin{document}


\begin{abstract}
  In this paper we describe the compactification of the Drinfeld
  modular curve.  This compactification is analogous to the
  compactification of the classical modular curve given by Katz and
  Mazur. We show how the Weil pairing on Drinfeld modules that we
  defined in earlier work gives rise to a map on the Drinfeld modular
  curve. We introduce the Tate-Drinfeld module and show how this
  describes the formal neighbourhood of the scheme of cusps of the
  Drinfeld modular curve.
\end{abstract}

\maketitle

\section{Introduction} \label{sec_intro}

\noindent 
Consider a smooth, projective,
geometrically irreducible curve $X$ over $\F_q$, and fix some point
$\infty$ on this curve. Let $$A:= \Gamma(X - \infty, \O_X)$$ be
the ring of functions on $X$ which are regular outside
$\infty$. Let $f \in A\bs \F_q$ be a non-constant element, and let 
$$A_f := A[f^{-1}].$$ The moduli
scheme $M^2(f)$ plays an important role in this paper. It represents
the functor which associates to every $A_f$-scheme $S$ the set of 
isomorphy classes of Drinfeld modules with a level $f$-structure over 
$S$.\\ In this paper we will address the following problems:

\begin{itemize}
  \item[$(i)$] Construct a morphism $$w_f: M^r(f) \lr M^1(f).$$ This morphism
    is induced by the Weil pairing for Drinfeld modules. The Weil
    pairing is defined in \cite{Hei02}. 
  \item[$(ii)$] Define the Tate-Drinfeld module and describe its universal
    property, using ideas of G. B\"ockle in Chapter $2$ of \cite{Boe02}
    and ideas of M. van der Put and J. Top in \cite{PT1} and \cite{PT2}.
  \item[$(iii)$] Describe a compactification $\ov M^2(f)$ of $M^2(f)$. 
    The Tate-Drinfeld module will enable us to describe the scheme of cusps
    $${\it Cusps} = (\ov M^2(f) - M^2(f))^{\rm red}.$$
  \item[$(iv)$] Compute the number of components of $M^2(f)$ 
    and describe the cusps of the analogue of the classical curve $X_0(N)$.
\end{itemize}

\noindent
We would like to point out that another description of the compactification 
of $M^2(f)$ can be found in T. Lehmkuhl's `Habilitation' \cite{Leh00}. 
Our treatment uses the Weil pairing and gives an explicit description of
the Tate-Drinfeld module. In a forthcoming
paper this enables us to develop the N\'eron model of the Tate-Drinfeld 
module, analogous to the Deligne and Rapoport's construction of the N\'eron 
model of the Tate-elliptic curve in \cite{DR73}. This will probably give
rise to an extension of the functor represented by $M^r(f)$ to a functor 
represented by the compactification $\ov M^r(f)$. 

\par\bigskip\noindent
Sections \ref{sec_defs} and \ref{sec_moduli} give a brief introduction
to the moduli problem and to the moduli schemes. At the end of Section
\ref{sec_moduli} we state the assumptions which are used throughout
this paper. In Section \ref{sec_defs} we recall the definition of
Drinfeld modules over schemes and level structures. In Section
\ref{sec_moduli} we describe the moduli problem that Drinfeld
considers in his original paper \cite{Drin74}.  The goal of Section
\ref{sec_weil-map} is to prove Theorem \ref{thm_map}, i.e., to
construct the morphism $w_f$ considered in problem $(i)$.  In Sections
\ref{sec_redth}, \ref{sec_TD} and \ref{sec_UP} we discuss problem
$(ii)$. In Section \ref{sec_redth} we classify the Drinfeld modules of
rank $2$ with level $f$-structure over the quotient field of some
complete discrete valuation ring $V$ which have stable reduction of
rank $1$. The main result of this section is Theorem \ref{thm_reduct}.
In Section \ref{sec_TD} we define the Tate-Drinfeld module of type
$\mf m$ with level $f$-structure. In Section \ref{sec_UP} we define
the universal Tate-Drinfeld module $\Cal Z$; cf.  Theorem
\ref{thm_STD_up}. For the application to problem $(iii)$ we need a
weak version of the universal property of the scheme $\Cal Z$ as
stated in Theorem \ref{thm_univTD}. In Sections \ref{sec_compact} and
\ref{sec_Cusp_TD} we study problem $(iii)$. In Section
\ref{sec_compact} we give a compactification of $M^2(f)$. The
treatment given here is analogously to the compactification of the
classical modular curve given by N.M. Katz and B. Mazur in Chapter 8
of their book \cite{KM85}. This enables us in Section
\ref{sec_Cusp_TD} to identify the formal neighbourhood of the scheme
of cusps $\ov M^2(f) - M^2(f)$ with the universal Tate-Drinfeld
module. The main results are Theorem \ref{thm_CompToTD} and
Corollary \ref{cor_cusp}. In Section \ref{sec_compon} we compute, using
the scheme of cusps, the number of components for every geometric fibre; 
cf. Theorem \ref{thm_geom_comp}.


\section{Drinfeld modules over schemes}
\label{sec_defs}

\noindent
Throughout this paper we will denote the quotient field of any integral
domain $D$ by $K_D$. We recall the definition of Drinfeld modules over 
schemes and level structures. There are many texts available for a more 
extensive account of these definitions; cf. \cite{Drin74}, \cite{Del87}, 
\cite{Sai96}, \cite{Leh00} and \cite{Tae02}.

\subsection{Line bundles and morphisms}
\noindent
Let $B$ be a commutative $\F_q$-algebra with $1$ and let $\G_{a, B}$
denote the additive group over $B$. The ring of 
$\F_q$-linear endomorphisms $\End_{\F_q}(\G_{a, B})$ of $\G_{a, B}$ is
isomorphic to the skew polynomial ring $B\{\t\}$. In this skew polynomial
ring multiplication is determined by the rule $\t b = b^q \t$ for all 
$b\in B$.

\par\bigskip\noindent
This can be generalized to schemes. Let $S$ be an $\F_q$-scheme, and let 
$L \lr S$ be a line bundle. As usual, $L$ is also a group scheme 
due to its additive group scheme structure. 
A {\em trivialization of $L$} is a covering $\Spec(B_i)$ of open affines 
of $S$ together with isomorphisms $L_{\mid \Spec(B_i)} \cong
\G_{a, B_i}$. By $\End_{\F_q}(L)$ we denote the
$\F_q$-linear $S$-group scheme endomorphisms of $L$. Let $\Cal L$ be
the invertible $\O_S$-sheaf corresponding to $L$, and let
$$\t^i: \Cal L \lr \Cal L^{q^i} \quad \mbox{by} \quad s
\mapsto s\otimes \ldots \otimes s.$$ The ring $\End_{\F_q}(L)$ is
isomorphic to the ring of all formal expressions $\sum_i \alpha_i \t^i$
which are locally finite. Here $\alpha_i: \Cal L^{q^i} \lr \Cal L$ is
an $\O_S$-module homomorphism for every $i$. Multiplication in the ring of
formal expressions is given by $\alpha_i\t^i \beta_j \t^j = \alpha_i \otimes
\beta_j^{q^i} \t^{i+j}$. If $\{ \Spec(B_i) \}_{i\in I}$ is a 
trivialization of $L$, then the restriction of $\End_{\F_q}(L)$ to 
$B_i$ is simply $B_i\{ \t \}$.\\
Furthermore, we denote by $\d_0$ the point derivation at $0$:
$$\d_0: \End_{\F_q}(L) \lr \Gamma (S, \Cal O_S) \quad \mbox{by} \quad
\sum_i \alpha_i \t^i \mapsto \alpha_0.$$

\subsection{Drinfeld modules over a scheme}

\begin{definition} \label{def_DMField}
  Let $K$ be an $A$-field equipped with an $A$-algebra structure
  given by $\g: A\lr K$. Let $L = \G_{a, K}$.
  A {\em Drinfeld module over $K$} is a ring homomorphism $$\ph:
  A \lr \End_{\F_q}(\G_{a, K})$$ such that 
  \begin{itemize}
    \item[$(i)$] $\d_0 \circ \ph = \g$;
    \item[$(ii)$] there is an $a \in A$ such that $\ph_a \neq \g(a)$.
  \end{itemize}
\end{definition}

\noindent
A Drinfeld module over a field $K$ has a {\it rank}, i.e., there is an
integer $r > 1$ such that $\deg_{\t} \ph_a = r \deg(a)$ for all $a
\in A$.

\begin{definition}
  Let $S$ be an $A$-scheme via the morphism $$\g_S: S \lr \Spec(A).$$ 
  A {\em Drinfeld module of rank $r$ over $S$} is a pair $(L, \ph)_S$
  of a line bundle $L\lr S$ and a ring homomorphism $\ph: A \lr
  \End_{\F_q}(L)$ such that
  \begin{itemize}
    \item[$(i)$] $\d_0 \circ \ph = \g_S^{\#}$;
    \item[$(ii)$] For all $a \in A$ the morphism $\ph_a$ is finite
      of degree $q^{r \deg(a)}$.
  \end{itemize}
\end{definition}

\begin{remark}
  The pull-back of a Drinfeld module $(L, \ph, S)$ along a morphism
  $$\Spec(K) \lr S$$ for some field $K$ is a Drinfeld module over
  $K$ in the sense of Definition \ref{def_DMField}.
\end{remark}

\noindent
If $S = \Spec(B)$ and $L$ is isomorphic to $\G_{a, B}$, then we simply
write $\ph$ instead of $(L, \ph)_S$. The morphism $\g_S$ is called the
{\it characteristic} of $(L, \ph)_S$. An ideal $\mf n \subset A$ is called 
{\it away from the characteristic} if $V(\mf n)$ is disjoint with the
image of $\g_S$. 


\begin{definition}
  A {\em morphism $\xi$ of Drinfeld modules over $S$} $$\xi: (L,
  \ph)_S \lr (M, \psi)_S$$ is a map $\xi \in \Hom_{\F_q}(L, M)$ such
  that $\xi \circ \ph_a = \psi_a \circ \xi$ for all $a \in A$.  A
  morphism $\xi$ is called an {\it isomorphism}, if it gives an isomorphism
  between the line bundles $L$ and $M$ over $S$.\\
  An {\it isogeny} of Drinfeld modules is a finite morphism.
\end{definition}

\begin{remark}
  An isogeny exists only between Drinfeld modules of the same rank. 
\end{remark}

\begin{remark}
  A Drinfeld module $(L, \ph)_S$ of rank $r$ and a morphism $$f: T
  \lr S$$ give by pull-back rise to a Drinfeld module $(f^* L, f^*
  \ph)_T$ over $T$ of rank $r$.
\end{remark}

\subsection{Level structures}
\noindent
For any non-zero $f \in A$, let $$\ph [f] := \ker(\ph_f: L \lr L).$$ This is a
finite, flat group scheme over $S$, namely $$\ker(\ph_f: L \lr L) = L
\times_L S$$ where the fibre product is taken over $\ph_f: L \lr L$ and 
the unit-section $e: S \lr L$ of the group scheme $L$. If $\mf n \subset
A$ is any non-zero ideal, then $$\ph[\mf n] = \prod_{f \in \mf n} \ph[f]$$ 
where the product is the fibre product over $L$. The scheme $\ph[\mf n]$ is
\'etale over $S$ if and only if $\mf n$ is away from the characteristic.

\begin{remark}
  Let $K$ be an algebraically closed $A$-field. 
  If $S = \Spec(K)$, then 
  $$\ph[\mf n] = \Spec(K[X]/(\ph_{f}(X), \ph_g(X)))$$ with $\mf n =
  (f, g)$.
\end{remark}

\begin{definition} \label{def_level}
  Suppose that $\mf n$ is away from the characteristic, then a level 
  $\mf n$-structure $\lambda$ over $S$ of $(L, \ph)_S$
  is an $A$-isomorphism $$\lambda : (A/\mf n)^r
  \stackrel{\sim}{\lr} \ph[\mf n](S).$$
\end{definition}

\begin{remark}
  If $\mf n$ is not away from the characteristic, then 
  the definition of a level $\mf n$-structure is more involved. 
  One can view $\ph[\mf n]$ as an
  effective Cartier divisor on $L$. By definition a level 
  $\mf n$-structure of $(L, \ph)_S$ is an $A$-homomorphism $$\lambda: (A/\mf
  n)^r \lr L(S)$$ which induces an equality of Cartier-divisors:
  $$\sum_{\alpha \in (A/\mf n)^r} \lambda(\alpha) = \ph[\mf n].$$ We
  will not expand on this. In this paper we restrict to the cases for
  which Definition \ref{def_level} is enough.
\end{remark}

\noindent
Let the triple $(L, \ph, \lambda)_S$ denote a Drinfeld module of rank
$r$ over $S$ with level $\mf n$-structure $\lambda$.

\begin{definition}
  A {\em morphism} between two triples $(L, \ph, \lambda)_S$ and $(M,
  \psi, \mu)_S$ is a morphism $\xi: (L, \ph)_S \lr (M, \psi)_S$ of
  Drinfeld modules over $S$ such that $\xi(S) \circ \lambda = \mu$ 
  where $\xi(S): L(S) \lr M(S)$ is induced by $\xi$. A
  morphism is called an {\em isomorphism} if $\xi$ is an isomorphism
  of Drinfeld modules.
\end{definition}



\section{The moduli problem}
\label{sec_moduli}

\noindent
Let $\mf n\subset A$ be a non-zero, proper ideal. 
In his original paper, Drinfeld considers the following moduli problem
for Drinfeld modules. Let $$\Cal F^r(\mf n): A-\mbox{\sc Schemes} \lr 
\mbox{\sc Sets}$$ be the functor which associates to each $A$-scheme $S$ 
the set of isomorphy classes of Drinfeld modules over $S$ of rank $r$ with
level $\mf n$-structure over $S$. Drinfeld showed the following;
cf. Proposition 5.3 and its corollary in \cite{Drin74}.

\begin{theorem}[V.G. Drinfeld] \label{thm_drin}
  If $\mf n \subset A$ is an ideal divisible by at least
  $2$ distinct primes, then there exists a fine moduli space $$M^r(\mf
  n) \lr {\rm Spec}(A)$$
  representing the moduli problem $\Cal F^r(\mf n)$. Moreover, this
  scheme has the following properties:
  \begin{itemize}
    \item[$(i)$] $M^r(\mf n)$ is affine and smooth of dimension $r$;
    \item[$(ii)$] $M^r(\mf n) \lr {\rm Spec}(A)$ is smooth of relative
      dimension $r-1$ over ${\rm Spec}(A) - V(\mf n)$.
  \end{itemize}
\end{theorem}

\noindent
For arbitrary non-zero ideals $\mf n \subset A$, the functor $\Cal 
F^r(\mf n)$ has in general only a coarse moduli scheme. This coarse
moduli scheme will also be denoted by $M^r(\mf n)$. We recall here 
briefly the construction of this scheme; cf. \cite{Tae02} for a nice 
exposition of this. Let $\mf n, \mf m \subset A$ be ideals such that $\mf n 
\mf m$ is divisible by at least two distinct prime ideals, then by Theorem
\ref{thm_drin}
there exists an affine scheme $M^r(\mf n \mf m)$ representing the moduli
functor $\Cal F^r(\mf n \mf m)$. Let $(L, \ph, \lambda)_S$ be a Drinfeld
module of rank $r$ with full level $\mf m \mf n$-structure $\lambda$ over
an $A$-scheme $S$. On these triples the group $\Gl_r(A/\mf m \mf n)$
acts by $$\s (L, \ph, \lambda) :=
(L, \ph, \lambda \circ \s) \quad \mbox{for all} \quad \s \in
\Gl_r(A/\mf n).$$
Note that $\Gl_r(A/\mf n)$ is isomorphic to the kernel of the 
$~{\rm mod}~\mf m$ reduction map $$\Gl_r(\mf m \mf n) \lr \Gl_r(\mf m).$$
This induces an action of $\Gl_r(A/\mf n)$ on $M^r(\mf m \mf n)$. The
coarse moduli scheme of $\Cal F^r(\mf m)$ is defined as $$M^r(\mf m)
:= M^r(\mf m \mf n) / \Gl_r(A/\mf n).$$ This quotient exists,
because $\Gl_r(A/\mf n)$ is finite. It is, however, not obvious that
this scheme is coarse for the given moduli problem. See \cite{Tae02}
for a proof of the coarseness of the scheme.
The scheme $M^r(\mf m)$ does not depend on the choice of $\mf n$. 

\subsection{Actions on $M^r(\mf n)$}
\noindent
Let $\hat A = \underset{\ll}{\lim}~A/\mf n$, and let $\A_f = \hat A 
\otimes_A K_A$. In this subsection we describe the natural action of 
$\A_f \cdot 
\Gl_r(\hat A)$ on $M^r(\mf n)$; cf. \cite[5D]{Drin74}. Using this action, 
we can define the action of $\Gl_r(A/\mf n)$ and ${\rm Cl}(A)$ on 
$M^r(\mf n)$. 
To keep this paper self-contained, we recall here the treatment given in 
Section 3.5 of \cite{Leh00}, where the reader can find proofs and details.

\par\bigskip\noindent
A {\it total level structure} of a Drinfeld module $(L, \ph)_S$ is a 
homomorphism $$\kappa: (K_A/A)^r \lr L(S)$$ such that its restriction to 
$(\mf n^{-1} A/A)^r$ defines a level $\mf n$-structure.
Let $$M^r := \underset{\ll}{\lim} ~M^r(\mf n)$$ where $\mf n$ runs through 
the non-zero ideals of $A$. This is an affine scheme, and $M^r$ represents
the functor which associates to each $A$-scheme $S$ the set of isomorphy
classes of Drinfeld modules with a total level structure over $S$. 

\par\bigskip\noindent
There is a natural action of $\Gl_r(\A_f)$ on $M^r$, which is defined as 
follows. Let $S$ be an $A$-scheme, and let $(L, \ph, \kappa)_S$ be a 
Drinfeld module with a total level structure over $S$. 
Let $\s \in \Gl_r(\A_f)$ such that the entries 
of $\s$ are elements of $\hat A$, then $\s$ gives rise to a map 
$$\s: (K_A/A)^r \lr (K_A/A)^r.$$ Let $H_{\s}$ denote the kernel of $\s$. 
The kernel $H_{\s}$ gives rise to a finite subgroup scheme 
of $L$. We can divide out the pair $(L, \ph)_S$ by this subgroup scheme.
This gives us a pair $(L', \ph')_S$. The following diagram equips
the pair $(L', \ph')_S$ with a total level structure $\kappa'$:

\begin{equation} \label{diag_ClGpAct}
  \begin{CD}
    0 @>>> H_{\s} @>>> (K_A/A)^r @> \s >> (K_A/A)^r @>>> 0\\
    & & @VVV @V \kappa VV @VV \kappa' V\\
    0 @>>> H_{\s}(S) @>>> L(S) @>>> L'(S).
  \end{CD}
\end{equation}

\noindent
If $\s$ comes from an element in $A\bs \{ 0 \}$, then its action is 
trivial. This implies that we get an action of $\Gl_r(\A_f)/K_A^*$
on $M^r(S)$. As this action is functorial in $S$, this defines 
an action of $\Gl_r(\A_f)/K_A^*$ on $M^r$. 

\par\bigskip\noindent
For the moduli scheme $M^r(\mf n)$ we have $M^r(\mf n) = \Gamma(\mf n) 
\bs M^r$ with $$\Gamma(\mf n) := \ker(\Gl_r(\hat A) \lr \Gl_r(A/\mf n)).$$
The restriction of the universal triple $(L, \ph, \kappa)$ on $M^r$
to $M^r(\mf n)$ gives the universal pair $(\ph, \lambda)$ on $M^r(\mf n)$
(Recall that the line bundle of the universal Drinfeld module on $M^r(\mf n)$ 
is trivial.)  
As $\A_f^* \cdot \Gl_r(\hat A)$ commutes with $\Gamma(\mf n)$ in 
$\Gl_r(\A_f)$, it 
follows that the action of $\Gl_r(\A_f)$ on $M^r$ defines an 
action of $\A_f^* \cdot \Gl_r(\hat A)$ on $M^r(\mf n)$. 
The normal subgroup $K_A^* \cdot \Gamma(\mf n) \subset 
\A_f^* \cdot \Gl_r(\hat A)$ acts trivially on $M^r(\mf n)$. 
Let $$G := \A_f^* \cdot \Gl_r(\hat A)/K_A^* \cdot \Gamma(\mf n).$$ 
As $\A_f^*/K_A^* \cdot \hat A^* \cong {\rm Cl}(A)$, it is not 
difficult to see that we have the following exact sequence
$$0 \lr \Gl_r(A/\mf n)/\F_q^* \lr G \lr {\rm Cl}(A) \lr 0.$$ 

\par\bigskip\noindent
To describe the action of $\Gl_r(A/\mf n)/\F_q^*$, let $\s \in 
\Gl_r(\hat A)$ and let $\tilde \s$ be the image of $\s$ under 
the reduction map $\Gl_r(\hat A) \lr \Gl_r(A/\mf n)/\F_q^*$. 
Then $$\s: (\ph, \lambda) \mapsto (\ph, \lambda \circ \tilde \s^{-1}).$$

\begin{remark}
  However, in the sequel we prefer to drop the inverse. If we talk about 
  the action of $\s \in \Gl_r(A/\mf n)/\F_q^*$ on $(\ph, \lambda)$, then we 
  mean the action given by 
  $$\s: (\ph, \lambda) \mapsto (\ph, \lambda \circ \s).$$ 
  Consequently, $\Gl_r(A/\mf n)/\F_q^*$ acts on the right of 
  $M^r(\mf n)$ and not on the left. 
\end{remark}

\par\bigskip\noindent
Let $m \in \hat A \cap \A_f^*$, then $m$ defines a unique ideal 
$\mf m = (m) \cap A \subset A$. We suppose that $\mf m$ is a non-zero, 
proper ideal which is relatively prime to $\mf n$.
Let ${\rm I}_r$ denote the identity element in 
$\Gl_r(\A_f)$, and let $\s = m \cdot {\rm I}_r$. We describe the
action of $\s$ on $M^r$. Clearly, $H_{\s} = (\mf m^{-1} A/A)^r$, and 
$H_{\s}$ maps to $\ph[\mf m](M^r)$ under $\kappa$. This means that the isogeny
$$\xi_{\mf m}: (L, \ph) \lr (L', \ph')$$ defined by $\s$ has kernel
$\ph[\mf m]$. The total level structure $\kappa'$ is given by
$$\kappa' = \xi_{\mf m} \circ \kappa \circ m^{-1}.$$
Let $\ph'$ denote the restriction of $(L', \ph')$ to $M^r(\mf n)$. 
Let $\ov m$ denote the image of $m$ under the reduction map $\hat A
\lr A/\mf n$. As $\mf m + \mf n = A$, we see that $\ov m \in (A/\mf
n)^*$. Let $\ph$ be the restriction of $(L', \ph')$ to $M^r(\mf n)$, 
then the action of $m$ on the universal pair $(\ph, \lambda)$ on 
$M^r(\mf n)$ is given by $$m: (\ph, \lambda) \mapsto (\ph', 
\xi_{\mf m} \circ \lambda \circ \ov m^{-1} ).$$
This describes the action of $m$ on $M^r(\mf n)$. 

\par\bigskip\noindent 
Let $\mf m \subset A$ be a non-zero ideal relatively prime to 
$\mf n$, i.e., $\mf m + \mf n = A$. Choose $m \in \hat A$ such that 
$(m) = \mf m \hat A$ and $m \equiv 1 ~{\rm mod}~ \mf n$. 
We define the action of $\mf m$ on 
$(\ph, \lambda)$ to be the action of $(m)$: 
$$\mf m: (\ph, \lambda) \mapsto (\ph', \xi_{\mf m} \circ \lambda).$$
This action is well-defined: the chosen element $m$ is unique up to an
element $\alpha \in \hat A^*$ with $\alpha \equiv 1 ~{\rm mod}~ \mf
n$. For such an element $\alpha$ we have $\alpha {\rm I}_r \in 
\Gamma(\mf n)$. Consequently, $m \cdot {\rm I}_r$ and $\alpha m \cdot 
{\rm I}_r$ give the same element in $G$.\\
Using this, we can define the {\it action of ${\rm Cl}(A)$ on $M^r(\mf n)$}.
First note that by Lemma \ref{lem_choice-m} we can represent every class 
in ${\rm Cl}(A)$ by a non-zero ideal $\mf m$ with $\mf m + \mf n = A$. 
Namely, suppose that $\mf m$ and $\mf m'$ are both non-zero
ideals relatively prime to $\mf n$ which represent the same class in
${\rm Cl}(A)$, then there is an element $x \in K_A^*$ with $\mf m = x
\mf m'$. Let $m, m' \in \hat A$ be elements which define the action of 
$\mf m$ and $\mf m'$, respectively. Then there is an element 
$\alpha \in \hat A^*$ with $\alpha \equiv 1 ~{\rm mod}~\mf n$ such that 
$m' = x\alpha m$ with $x \alpha \cdot {\rm I}_r \in K_A^* \cdot 
\Gamma(\mf n)$. Therefore, $(m)$ and $(m')$ give the same element in $G$. 

\begin{remark}
  These considerations also imply that 
  $$G\cong {\rm Cl}(A) \times \Gl_r(A/\mf n)/\F_q^*.$$
\end{remark}

\begin{remark}
  Let us once more stretch that we use the convention that the action 
  of $\s \in \Gl_r(A/\mf n)$ on 
  $M^r(\mf n)$ is given by $$\s: (\ph, \lambda) \lr (\ph, \lambda \circ \s).$$ 
  So $G = {\rm Cl}(A)\times \Gl_r(A/\mf n)/\F_q^*$ acts on $M^r(\mf n)$
  on the right.
\end{remark}

\subsection{Assumptions in this paper}

\noindent
In this paper we will make the following two assumptions.

\begin{enumerate}
  \item Throughout this paper we will assume that $\mf n = (f)$ is a 
    non-zero, proper, principal ideal. This simplifies the description of the
    Tate-Drinfeld module. Dropping this assumption does not seem to 
    give rise to different results.\\ 
    If $\mf m\subset A$ is a non-zero
    proper ideal containing $f$, then by the previous, we see that
    $M^r(\mf n) = M^r(f)/G$ where $G$ is given by dividing out the action 
    of the kernel $$\ker(\Gl_r(A/fA) \lr \Gl_r(A/\mf n)).$$ 
  \item We will not consider the moduli problem over $A$-schemes, but
    over $A_f$-schemes $S$, i.e., $f$ is invertible in $S$. This implies
    that $f$ is away from the characteristic of $(L, \ph)_S$.
    This assumption is used because the Weil pairing plays an important
    role in our description, and the Weil pairing is in \cite{Hei02} only
    defined for $f$-torsion which is away from the 
    characteristic.\\ This assumption implies that a level $f$-structure 
    is an isomorphism $$\lambda: (A/fA)^r \stackrel{\sim}{\lr} \ph[f](S).$$
\end{enumerate}

\noindent
So in this paper we will be considering the moduli problem
$$\Cal F^r(f): A_f-\mbox{\sc Schemes} \lr \mbox{\sc Sets}$$ which
associates to each $A_f$-scheme $S$, the set of isomorphy classes of
Drinfeld modules of rank $r$ with full level $f$-structure
over $S$. We write $M^r(f)$ for the scheme which represents
$\Cal F^r(f)$. It follows from the proof of Theorem \ref{thm_drin} 
that the moduli scheme $M^r(f)$ is a fine moduli scheme if $f \neq 1$.

\par\bigskip\noindent
Throughout this paper we will write $\Spec(R) = M^1(f)$. The ring $R$
is regular and $M^1(f)$ is connected. In fact, $R$ is the integral 
closure of $A_f$ in a field extension of $K_A$. The Galois group of 
$K_R/K_A$ is the group $$G \cong (A/fA)^*/\F_q^* \times {\rm Cl}(A)$$ 
that we discussed above; cf. Section 8 in \cite{Drin74}.


\section{The Weil pairing on the modular schemes.} \label{sec_weil-map}

\noindent 
In this section we will show that the Weil pairing for
Drinfeld modules over an $A_f$-field $K$ as defined in the
previous paper gives rise to the following theorem:

\begin{theorem} \label{thm_map}
  The Weil pairing induces an $A_f$-morphism $$w_f: M^r(f)
  \lr M^1(f).$$ The Weil pairing is ${\rm Cl}(A) \times 
  \Gl_r(A/fA)$-equivariant. 
\end{theorem}

\noindent
Let $(\ph, \lambda)$ be a Drinfeld module $\ph$ of rank $r$ over $K$ with
level $f$-structure. The Weil pairing is an $A/fA$-isomorphism $$w_f:
\wedge^r \ph[f](K) \stackrel{\sim}{\lr} \psi[f](K) \otimes_A
\Omega_f^{\otimes r-1}.$$ It is unique up to a unique isomorphism of 
$\psi$. Once and
for all we fix a generator $\omega$ of the $A/fA$-module
$\Omega_f$. This gives an $A/fA$-isomorphism $$w_f: \wedge^r
\ph[f](K) \stackrel{\sim}{\lr} \psi[f](K).$$

\par\smallskip\noindent
The level $f$-structure $\lambda$ induces a canonical isomorphism
$$\wedge^r \lambda: \wedge^r (A/fA)^r \stackrel{\sim}{\lr}
\wedge^r \ph[f](K).$$ Because $\wedge^r (A/fA)^r$ is canonically
isomorphic to $A/fA$, $\psi$ comes equipped with a level $f$-structure
$\mu$ over $K$ via the following commutative diagram:
$$\begin{CD}
  \wedge^r (A/fA)^r @> \wedge^r \lambda >> \wedge^r \ph[f](K)\\
  @AAA @V w_f VV\\
  A/fA @> \mu >> \psi[f](K).
\end{CD}$$

\noindent
Note that if $\xi$ is an isomorphism between $\psi$ and $\psi'$, then
also the pairs $(\psi, \mu)$ and $(\psi', \mu')$ are isomorphic via
$\xi$. Here $\mu$ and $\mu'$ are defined by the previous diagram by
$\psi$ and $\psi'$ respectively. So the pair $(\psi, \mu)$
is unique up to isomorphy. These considerations show
the following: 

\begin{lemma} \label{lem_isolev}
  The Weil pairing gives for all $A_f$-fields $K$ rise to a map
  $$w_K: M^r(f)(K) \lr M^1(f)(K).$$
  This map depends functorially on $K$.
\end{lemma}
\begin{proof}
  The construction of the map $w_K$ is described above. It
  associates to each isomorphy class $(\ph, \lambda)$ of rank $r$ over
  $K$ a unique isomorphy class $(\psi, \mu)$ of rank $1$ over $K$.
  That this map depends functorially on $K$ follows immediately from
  the construction in terms of $A$-motives.
\end{proof}

\noindent
From this lemma, we proceed as follows to prove the existence of the map
$w_f$. 

\begin{proof}[Proof of Theorem \ref{thm_map}] 
  Let $(\ph, \lambda)$ be the universal pair over $M^r(f)$.
  The moduli scheme $M^r(f)$ is affine and regular over
  $\Spec(A_f)$. So we may write $$M^r(f) \cong \coprod_{i=0}^n
  \Spec(S_i)$$ such that
  each $S_i$ is an integrally closed domain of relative dimension $r-1$
  over $A_f$. Moreover, from Drinfeld's description we know that
  $$M^1(f) = \Spec(R)$$ where $R$ is an integrally closed domain. 
  Let $K_{S_i}$ be the quotient field of
  $S_i$. Lemma \ref{lem_isolev} gives rise to a unique
  isomorphy class $$(\psi, \mu) \in M^1(f)(\prod_i K_{S_i}).$$ This means
  that there exists a unique $A_f$-ring homomorphism $$h: R
  \lr \prod_j K_{S_j}.$$ Hence, the theorem 
  follows if we can show that $h(R) \subset 
  \prod_j S_j$.\\ Let $\mf p \in \Spec(S_j)$
  be a closed point of height $1$. By Lemma \ref{lem_reg} it follows
  that there is a map $\Spec(S_{j, \mf p}) \lr M^1(f)$, inducing
  the map $$\Spec(K_{S_j}) \lr M^1(f).$$ Consequently, 
  $$h(R) \cap K_{S_j} \subset S_{j, \mf p}.$$ As 
  $S_j = \cap_{\mf p} S_{j, \mf p}$
  where the intersection runs over all primes of height $1$ in $S_j$, one has
  $h(R) \cap K_{S_j} \subset S_j$.\\
  For the ${\rm Cl}(A) \times \Gl_r(A/fA)$-equivariance of $w_f$, note
  that the $\Gl_r(A/fA)$-equivariance is obvious. Let $\mf m \subset 
  A$ be a non-zero, proper ideal with $f \not \in \mf m$ representing 
  an element in ${\rm Cl}(A)$. For the 
  ${\rm Cl}(A)$-equivariance we recall its definition in the 
  previous section. Using the notations
  of the previous section, we see that the action of $\mf m$
  on $(\ph, \lambda)$ is given by $$(\ph, \lambda) \mapsto (\ph', 
  \xi_{\mf m} \circ \lambda),$$ with $\xi_{\mf m} \ph_a = \ph'_a \xi_{\mf m}$
  for all $a \in A$.   
  Let $(\psi, \mu)$ be the image of $(\ph, \lambda)$ under the 
  Weil pairing, and let $(\psi', \mu')$ be the image of $(\ph', 
  \xi_{\mf m} \circ \lambda)$ under $w_f$.
  Let $F = \prod_i \ov K_{S_i}$. The isogeny $\xi_{\mf m}$ induces an 
  isogeny $\zeta_{\mf m}: \psi \lr \psi'$. The kernel of $\zeta_{\mf m}$
  is  
  $$\ker(\zeta_{\mf m})(F) = \wedge^r \ker(\xi_{\mf m})(F) 
  = \wedge^r \ph[\mf m](F) \cong \psi[\mf m](F).$$ Therefore, the 
  action of $\zeta_{\mf m}$ on $\psi$ coincides with the action of 
  $\mf m$ on $\psi$. So we have $$\mf m: (\psi, \mu) \mapsto 
  (\psi', \zeta_{\mf m} \circ \mu) = (\psi', \mu').$$
\end{proof}

\begin{lemma} \label{lem_reg}
  Let $S$ be a regular local $A_f$-ring, let $K_S$ be its quotient
  field, and let $(\ph, \lambda) \in M^r(f)(S)$. The unique class
  $(\psi, \mu) \in M^1(f)(K_S)$ associated to $(\ph, \lambda)$ by the
  Weil pairing comes from a unique class $(\psi', \mu') \in M^1(f)(S)$
  via the canonical embedding $S \lr K_S$.
\end{lemma}
\begin{proof}
  Via the ring homomorphism $S \lr K_S$ we may view the pair
  $(\ph, \lambda)$ over $K_S$, and we can associate via the Weil pairing
  a pair $(\psi, \mu)$ of rank $1$ over $K_S$. We want to prove
  that there is a representing pair in the isomorphy class of
  $(\psi, \mu)$ which is defined over $S$.\\
  Let $V$ be the set of all height one primes of $S$.
  The ring $S$ is a UFD; cf. Theorem 20.3 in \cite{Mat80}.
  Consequently, every $\mf p \in V$ is of the form $\mf p = (h_{\mf p})$
  for some irreducible $h_{\mf p} \in S$; cf. Theorem 20.1 in
  \cite{Mat80}. Let $v_{\mf p}$ denote the valuation at $\mf p$.
  The invertible elements of $S$ are given by $$S^* = \{ s \in S
  \mid v_{\mf p}(s) = 0 ~\mbox{for all $\mf p \in V$} \}.$$
  Let $S_{\mf p}$ denote the local ring of $S$ at $\mf p$. This is
  a discrete valuation ring.\\
  As the $f$-torsion of $\psi$ is $K_S$-rational, it follows that
  $\psi$ has good reduction at every $\mf p \in V$. This implies
  that for all $\mf p \in V$ there is an element $x_{\mf p}
  \in K_S$ such that $x_{\mf p} \psi x_{\mf p}^{-1}$ is a Drinfeld module
  of rank $1$ defined over $S_{\mf p}$. In fact, we may assume
  $x_{\mf p} = h_{\mf p}^{m_{\mf p}}$ for some $m_{\mf p} \in \Z$.\\
  There are only finitely many $m_{\mf p} \neq 0$, as we show below.
  So we can define $x = \prod_{\mf p \in V}
  h_{\mf p}^{m_{\mf p}}$. The pair $(x \psi x^{-1}, x \mu)$ is defined
  over $S$. So this is the pair that we are looking for.\\
  To see that there are only finitely many $m_{\mf p} \neq 0$,
  let $a \in A \bs \F_q$ and consider the leading coefficient $c$
  of $\psi_a$. Under the isomorphism $x_{\mf p}$, the leading coefficient
  becomes $x_{\mf p}^{1 - q^{\deg(a)}} c$. As $v_{\mf p}(c) = 0$ for
  all but finitely many $\mf p$'s, it follows that $m_{\mf p} = 0$ for
  all but finitely many $\mf p$'s.
\end{proof}

\begin{prop}
  The morphism $w_f$ induces the maps $w_K$ for any $A_f$-field $K$,
  where $w_K$ is as in Lemma \ref{lem_isolev}.
\end{prop}
\begin{proof}
  By the functoriality in $K$ of the maps $w_K$, it suffices to prove
  the statement for algebraically closed fields $\ov K$. So let $\ov K$ be
  algebraically closed and let $$\zeta : \Spec(\ov K) \lr M^r(f)$$ be a
  geometric point of $M^r(f)$. If $\zeta$ is a generic point of one of the
  connected components of $M^r(f)$, then $w_f$ induces $w_{\ov K}$ by
  construction. If $\zeta$ is a closed point we have to do something more.
  Clearly, $\zeta$ factors over $\Spec(k_{\zeta})$ where 
  $k_{\zeta}$ is the residue field at (the image of) $\zeta$. To see that
  $w_{k_{\zeta}}$ and consequently $w_{\ov K}$ is induced by $w_f$, we have
  to dive into the language of $A$-motives a little; cf. \cite{Hei02}.\\
  Let $V = \O_{\zeta}$, then $V$ is a regular local $A_f$-ring and
  let $(\ph, \lambda)$ be defined over $V$ of rank $r$. Let $K_V$ be the
  quotient field of $V$. Then the
  construction of $A$-motives associates to $\ph$ the Drinfeld module
  $\psi$ of rank $1$ for which $$M(\psi) \cong \wedge^r M(\ph).$$ By
  Lemma \ref{lem_reg}, the pair $(\psi, \mu)$ is also defined over
  $V$. Because $\ph$ is defined over $V$, it makes sense to consider
  $M^0(\ph) = V\{ \t \}$ with the obvious $V\{ \t \}\otimes_{\F_q} A$-action
  such that $$M^0(\ph) \otimes_V K_V = M(\ph).$$ Similarly, we can define
  $M^0(\psi)$ because $\psi$ is defined over $V$. Clearly,
  $$\wedge_{K_V \otimes A}^r M(\ph) \cong K_V \otimes_V \wedge_{V\otimes
  A}^r M^0(\ph).$$ Consequently, the $K_V\{\t\} \otimes A$-isomorphism
  $$M(\psi) \cong \wedge^r M(\ph)$$ comes from a $V\{ \t \}\otimes
  A$-isomorphism $$M^0(\psi) \cong \wedge_{V\otimes A} M^0(\ph).$$
  This construction can be reduced modulo the maximal of $V$. This
  gives us the construction of $w_{k_{\zeta}}$. Therefore, $w_f$
  induces $w_{\ov K}$.
\end{proof}


\section{Drinfeld modules of rank $2$ with stable reduction of rank $1$}
\label{sec_redth}

\noindent
For this section we fix the following notation. Let $V$ be a complete
discrete valuation ring which is also an $A_f$-algebra. Let $K_V$ be the
quotient field of $V$, and let $\pi \in V$ be a generator of the maximal
ideal of $V$. Let $v(x)$ denote the $\pi$-valuation of $x$ for every 
$x \in V$. 

\begin{definition}
  Let $\ph$ be a Drinfeld module of rank $r$ over $K_V$, then $\ph$
  has {\em stable reduction at $v$ of rank $r'$} if
  $\ph$ is isomorphic over $K_V$ to a Drinfeld module $\ph'$ over $K_V$
  such that for all $a \in A$ each coefficient $\beta_i(a)$ of the sum 
  $\ph'_a = \sum \beta_i(a) \t^i$ is an element of $V$ and the 
  reduction $\ph' ~{\rm mod}~\pi V$ is a Drinfeld module of rank $r'$ over
  $V/\pi V$.\\
  The Drinfeld module $\ph$ has {\em potentially 
  stable reduction at $v$ of rank $r'$} if there is a finite field 
  extension $L$ of $K_V$ and a valuation $w$ of $L$ extending $v$ such 
  that $\ph$ has stable reduction at $w$ of rank $r'$. 
\end{definition}

\noindent
Let $\ph$ be a Drinfeld module of rank $2$ over $K_V$
with full level $f$-structure $\lambda$ over $K_V$ such that $\ph$
has potentially stable reduction of rank $1$.
The goal of this section is Theorem \ref{thm_reduct} which describes
the pairs $(\ph, \lambda)$ in terms of Drinfeld modules of rank $1$
and lattices.

\begin{lemma} \label{lem_stab_red}
  Let $\ph$ be a Drinfeld module with $K_V$-rational $f$-torsion. 
  If $\ph$ has potentially stable reduction at $(\pi)$, then $\ph$
  has stable reduction at $(\pi)$. 
\end{lemma}
\begin{proof}
  As $\ph$ has potentially stable reduction, there is an element 
  $k \in \Q$ such that for all $x \in K_V^{\rm sep}$ with $v(x) = k$  
  we have the following: $\tilde \ph := x \ph x^{-1}$ is a Drinfeld
  module, $\tilde \ph_a$ has all coefficients in $V$ for all $a \in A$,
  and $\tilde \ph ~{\rm mod}~ \pi V$ is a Drinfeld module over $V/\pi V$. 
  Cf. Section 4.10 in \cite{Gos96}.
  As $f \in V^*$, it is not difficult to see that $k$ is the smallest
  slope of the Newton polygon of $\frac{1}{X} \ph_f(X)$. Therefore, 
  $$k = \max \{ v(\alpha) \mid \alpha \in \ph[f](\ov K_V)\bs \{ 0 \} \}.$$ 
  As the $f$-torsion of $\ph$ is $K_V$-rational, we have $k \in \Z$. 
  Therefore, we may choose $x \in K_V$. 
\end{proof}

\noindent
To abbreviate notation, we introduce the following two properties $P$
and $P'$. Let $\chi$ be a Drinfeld module of rank $r$ over $K_V$.

\begin{itemize}
  \item[$P(\chi)~:$] $\chi$ has stable reduction of rank $1$.
  \item[$P'(\chi):$] $\chi$ has stable reduction of rank $1$, $\chi_a$
    has all coefficients in $V$ for all $a \in A$ and $\chi[f](V) \cong A/fA$.
\end{itemize}

\noindent
Suppose that $\kappa$ is a level $f$-structure of $\chi$ over $K_V$.
If $r = 1$, then $P'(\chi)$ implies that $(\chi, \kappa)$ is defined
over $V$.

\par\bigskip\noindent
By Lemma \ref{lem_stab_red} we have $P(\ph)$ for the pair $(\ph, \lambda)$.
Therefore, we may assume that $\ph_a$ has all its coefficients in $V$
for all $a \in A$. Moreover, we have that the smallest slope
of the Newton polygon of $\ph_f$ equals $v(\ph_f) = 0$.
Consequently, $\ph[f](V) \cong A/fA$.
So in the isomorphy class of $(\ph, \lambda)$, there is a
representing element $(\ph, \lambda)$ with $P'(\ph)$ and this element
is unique up to $V^*$. In fact,
\begin{equation} \label{eq_iso1}
  \{ (\ph, \lambda)_{K_V} \mbox{with $P'(\ph)$} \} /V^*
  \stackrel{\rm bij}{\cong} \{ (\ph, \lambda)_{K_V} \mbox{with $P(\ph)$}
  \} / K_V^*.
\end{equation}

\par\bigskip\noindent
Note that the Weil pairing equips $V$ with 
an $R$-structure. The isomorphy class of $(\ph, \lambda)$ is induced by an 
$A_f$-morphism $\Spec(K_V) \lr M^2(f)$. Composing this morphism with 
$w_f$ gives rise to an $A_f$-linear ring homomorphism $R \lr K_V$. 
As $R$ is integral over $A_f$, it follows that this ring homomorphism 
gives an $A_f$-linear ring homomorphism $h: R \lr V$.

\subsection{Drinfeld's bijection without level structure}
\noindent
To classify the isomorphy classes $(\ph, \lambda)$ with stable reduction
of rank $1$, we recall Drinfeld's classification of
Drinfeld modules of rank $2$ with potentially stable
reduction of rank $1$; cf. Proposition 7.2 in \cite{Drin74}.\\
An {\it $A$-lattice of rank $1$ in $K_V^{\rm sep}$} is a
projective $A$-module of rank $1$ which lies discretely in $K_V^{\rm
sep}$ and which is invariant under the action of $G_{K_V}:=
\Gal(K_V^{\rm sep}/K_V)$. Two $A$-lattices $\Lambda_1$ and
$\Lambda_2$ are called {\it isomorphic} if there is an element 
$x \in (K_V^{\rm sep})^*$
such that $x \Lambda_1 = \Lambda_2$. Then Drinfeld's result states the
following:

\begin{theorem}[V.G. Drinfeld] \label{thm_Dri}
  There is a bijection between the set of isomorphy classes over $K_V$
  of Drinfeld modules of rank $2$ over $V$ with potentially stable
  reduction of rank $1$ and the set of isomorphy classes over $K_V$ of
  pairs $(\psi, \Lambda)$, where $\psi$ is a Drinfeld module of rank $1$
  over $V$ and $\Lambda$ is an $A$-lattice of rank $1$ inside
  $K_V^{\rm sep}$.
\end{theorem}
\begin{proof}[Sketch of the proof.]
  Applying Proposition 5.2 in \cite{Drin74} to the ring $V_n = V/(v)^n$
  with $n\in \Z_{\geq 1}$ gives us the existence of unique elements
  $s_n \in V_n\{\t\}$ such that $s_n \ph_a s_n^{-1}$ is in standard
  rank $1$ form over $V_n$. Moreover, each $s_n$ has the form $$s_n = 1
  + \sum_{i=1}^{k_n} v_i \t^i \quad \mbox{with} \quad v_i \in (v).$$
  Let $s = \underset{\ll}{\lim} ~s_n$, then $s$ is an element in
  $V\{\! \{\t\}\! \}$, the set of skew formal power series in $\t$ over $V$.
  In the proof of Proposition 7.2 in \cite{Drin74}, Drinfeld shows that the
  homomorphism $s$ is in fact analytic. One has by construction
  $$s = 1 + \underset{i \geq 1}{\sum} v_i \t^i \quad \mbox{with}
  \quad v_i \in (v).$$ This implies both that $\Lambda:= \ker(s)$ is
  contained in $K_V^{\rm sep}$. Moreover, each element in
  $\Lambda \backslash \{ 0 \}$ has strictly negative valuation.\\ Let
  $\psi' = s \ph s^{-1}$, then $\psi' \mod (v) = \psi$. 
  We get the following diagram, which is commutative for all $a \in A$:

  \begin{equation} \label{diag1}
  \begin{CD}
  0 @>>> \Lambda @>>> K_V^{\rm sep} @> e_{\Lambda} >> K_V^{\rm sep}\\
  & & @VV \psi_a V @VV \psi_a V @VV \ph_a V\\
  0 @>>> \Lambda @>>> K_V^{\rm sep} @> e_{\Lambda} >> K_V^{\rm sep},
  \end{CD}
  \end{equation}

  \noindent
  where
  $$e_{\Lambda}(z) = z \prod_{\alpha \in \Lambda \backslash \{0\}}
  \left( 1 - \frac{z}{\alpha} \right) = s(z).$$
  Let $a \in A \backslash \F_q$, then $\psi_a^{-1}\Lambda$ is mapped
  surjectively to $\ph[a]$, hence $\psi_a^{-1} \Lambda/\Lambda \cong
  (A/aA)^2$ as $A$-module. On the other hand, the kernel of the
  surjective map $$\psi_a^{-1} \Lambda/\Lambda \lr \Lambda/\psi_a
  \Lambda$$ is isomorphic to $\psi[a](K_V^{\rm sep})$. So we may
  conclude that $\Lambda/\psi_a \Lambda \cong A/aA$. This implies that
  $\Lambda$ is a projective $A$-module of rank $1$. We already saw that
  $\Lambda$ consists of elements with strictly negative valuation, hence
  $\Lambda$ lies discretely in $K_V^{\rm sep}$. Finally, for any element
  $\s \in G_{K_V}$, we have $\s \circ s(z) = s \circ \s(z)$, hence
  $\Lambda$ is $G_{K_V}$-invariant. We conclude that $\Lambda$ is an 
  $A$-lattice of rank $1$ in $K_V^{\rm sep}$.
\end{proof}

\subsection{The bijection with level $f$-structure}

\noindent
Let $(\psi^{\rm un}, \mu^{\rm un})$ be the universal pair of rank $1$ 
defined over $R$. We will assume that $\mu^{\rm un}(1) = 1$, which is
possible because $f$ is invertible in $R$.
Drinfeld's construction in Theorem \ref{thm_Dri} lifts the rank 
$1$ Drinfeld module 
$\ph~{\rm mod}~\pi V$ to a unique Drinfeld module $\psi$ of 
rank $1$ defined over $V$. Also, the $f$-torsion of $\psi$
is $V$-rational. We would like to equip $\psi$ with a natural 
level $f$-structure $\mu$ which comes from $\lambda$. For this
we use the ring homomorphism $$h: R \lr V$$ which arises from the 
Weil pairing.\\
Suppose that $\psi$ is equipped with a level $f$-structure $\tilde \mu$, 
then the isomorphy class of $(\psi, \tilde \mu)$ comes from an $A_f$-linear 
ring homomorphism 
$$\tilde h: R \lr V,$$ i.e., there is a unique element $v \in V^*$ such that 
$$\tilde h( (\psi^{\rm un}, \mu^{\rm un}) ) = (v \psi v^{-1}, v \tilde \mu).$$ 
The pair $(\tilde h, v)$ uniquely determines $(\psi, \tilde \mu)$.

\par\bigskip\noindent
As $R$ is integral over $A_f$, there exists an $A_f$-automorphism 
$$g_{\mu}: R \lr R$$ with $g_{\tilde \mu} \in G = \Gal(R/A_f)$ such 
that $\tilde h = h \circ g_{\tilde \mu}.$ We
described this Galois group in Section \ref{sec_moduli}: 
$$G \cong {\rm Cl}(A) \times (A/fA)^*/\F_q^*.$$ 
So the element $g_{\tilde \mu}$ is given by a pair 
$(\mf m, \s)\in {\rm Cl}(A) \times (A/fA)^*/\F_q^*$.\\
If we have another level $f$-structure $\mu'$, then $\mu' = \alpha \tilde \mu$
for some $\alpha \in (A/fA)^*$ and we see that $g_{\mu'}$ corresponds
to the pair $(\mf m, \alpha \s)$.

\par\bigskip\noindent
We equip $\psi$ with the level $f$-structure $\mu$ such that 
$g_{\mu}$ is given by the pair $(\mf m, 1)$. More specifically, let 
$h_{\mu} := h \circ g_{\mu}$. Let $v \in V^*$ be an element with 
$v \psi v^{-1} =  \tilde h(\psi^{\rm un})$, then this element $v$ is unique
up to $\F_q^*$. Let $(\psi, \mu)$ be the pair determined by $(\tilde h, v)$.\\ 
The pair $(\psi, \mu)$ that we obtain in this way is unique up to $\F_q^*$. 
In particular, the isomorphy class of $(\psi, \mu)$
is uniquely determined in this way.

\begin{remark} \label{rem_vis}
  To make the choice of an element in $\F_q^*$ in the above
  construction `visible', we consider the set of elements $z$ in
  $\ph[f](K_V)\bs \ph[f](V)$ with $$w_f(z, e_{\Lambda}(\mu(1))) \in
  \F_q^* \cdot w_f(\lambda(1,0), \lambda(0,1)).$$ The set of these elements
  equals $\F_q^* \cdot z'$ where $z'$ is any element in this set. 
  Fix one of those
  elements $z$. We can now equip $\psi$ with the unique $\mu$ such 
  that $$w_f(z, e_{\Lambda}(\mu(1))) = w_f(\lambda(1,0), \lambda(0,1)).$$
  In this way, we can associate to a triple $(\ph, \lambda, z)$ a unique
  triple $(\psi, \mu, \Lambda)$.
\end{remark}

\par\bigskip\noindent
Our next goal is to prove Theorem \ref{thm_reduct}.
Drinfeld's construction gives an $A$-lattice $\Lambda$ of rank $1$ 
such that the pair $(\psi, \Lambda)$ determines $\ph$ and vice versa.\\
As in Drinfeld's proof, we can consider $\psi_a$ as a map
$$\psi_a: K_V^{\rm sep} \lr K_V^{\rm sep}$$ for every $a \in A$.
Because $f \not \in \ker(A_f \lr K_V)$, it follows that all roots of
the equation $\psi_a(X) = \alpha$ lie in $K_V^{\rm sep}$ for all $\alpha\in
\Lambda$. And thus we have $$(\psi_f)^{-1} \Lambda \subset K_V^{\rm sep},$$ as
in Drinfeld's proof.

\begin{lemma} \label{lem_rat-lat}
  If the $f$-torsion of $\ph$ is $K_V$-rational, then $(\psi_f)^{-1} \Lambda
  \subset K_V$. In particular $\Lambda \subset K_V$.
\end{lemma}
\begin{proof}
  Cf. Lemma 2.4 in \cite{Boe02}. By definition $\Lambda$ is
  invariant under $G_{K_V}$. Consequently, also $(\psi_f)^{-1}
  \Lambda$ is invariant under $G_{K_V}$ as the coefficients of
  $\psi_f$ lie in $V$. In fact, the action of $G_{K_V}$ on
  $(\psi_f)^{-1}\Lambda$ splits according to the following splitting exact
  sequence of $A$-modules: $$0 \lr \psi[f](K_V) \lr (\psi_f)^{-1}
  \Lambda \lr ((\psi_f)^{-1} \Lambda)_{\rm proj} \lr 0.$$ The $A$-module
  $((\psi_f)^{-1} \Lambda)_{\rm proj}$ is the projective part of
  $(\psi_f)^{-1}\Lambda$, which is isomorphic to $\Lambda$.\\
  By Drinfeld's construction there is an analytic map $e_{\Lambda}$
  with the commuting property as in diagram (\ref{diag1}), and this map
  commutes with the action of $G_{K_V}$. Moreover, via this map
  $\ph[f]$ is isomorphic to $(\psi_f)^{-1}\Lambda/\Lambda$. By the
  assumption that $\ph[f]$ is $K_V$-rational, it follows that
  $G_{K_V}$ acts trivially on $(\psi_f)^{-1}\Lambda/\Lambda$ and thus on
  $\psi[f]$. The latter fact implies that the only action of
  $G_{K_V}$ on $(\psi_f)^{-1}\Lambda$ is on the projective part of
  $(\psi_f)^{-1}\Lambda$, hence this action gives a subgroup $G$ of $\Gl_1(A)
  = \F_q^*$. On the other hand, this subgroup $G$ maps injectively
  into $\Aut_A((\psi_f)^{-1}\Lambda/\Lambda)$. But $G_{K_V}$ acts
  trivially on $(\psi_f)^{-1}\Lambda/\Lambda$, hence $G$ is trivial.
\end{proof}

\noindent
So we see that we can associate to a pair $(\ph, \lambda)$ with 
$P'(\ph)$ a triple $(\psi, \mu, \Lambda)$ such that $(\psi, \mu)$ is defined 
over $V$ and $\Lambda$ is an $A$-lattice with $(\psi_f)^{-1}
\Lambda \subset K_V$. This triple is unique up to $\F_q^*$.
I.e., the pair $(\ph, \lambda)$ determines a unique element in  
$\{ (\psi, \Lambda, \mu) \} /\F_q^*$.

\par\bigskip\noindent
Note, however, that the level $f$-structure $\lambda$ is not the
unique level structure such that $(\ph, \lambda)$
is mapped to this unique element in 
$\{ (\psi, \mu, \Lambda \} /\F_q^*$. In fact, if $\s \in 
\Gl_2(A/fA)$, then $(\ph, \lambda \circ \s) \mapsto \{ (\psi, \mu,
\Lambda) \} /\F_q^*$ if and only if $\det(\s) \in \F_q^*$. 
This is due to the fact that the Weil pairing is $\Gl_2(A/fA)$-equivariant.
Therefore, if one changes $\lambda$ by $\s$, then the morphism $h$
induced by the Weil pairing $w_f$ changes by $\det(\s)$. Consequently, 
the triple $(\psi, \Lambda, \mu)$ changes by $\det(\s)$. 
Let $\Sigma$ denote the subgroup of $\Gl_2(A/fA)$ given by 
$$\Sigma = \{ \s \in \Gl_2(A/fA) \mid \det(\s) \in \F_q^* \}.$$
The previous shows that the map 
$$\{ (\ph, \lambda) \} / \Sigma \lr \{ (\psi, \mu, \Lambda) \} /\F_q^*$$
is injective. 

\begin{remark} \label{rem_injpair}
  Again, once we have chosen an element $z$ as in Remark \ref{rem_vis}, 
  we get an injective map 
  $$\{ (\ph, \lambda, z) \} /\Sl_2(A/fA) \lr  \{ (\psi, \mu, \Lambda) \}.$$
\end{remark}

\par\bigskip\noindent
On the other hand, let a triple $(\psi, \mu, \Lambda)$ be given. 
We show that there exists a pair $(\ph, \lambda)$ such that 
under the above construction $(\ph, \lambda)$ is mapped to 
the class $\{ (\psi, \mu, \Lambda) \}/\F_q^*$.\\ 
The triple $(\psi, \mu, \Lambda)$ gives rise to the following:

\begin{enumerate}
  \item a morphism $\tilde h: R \lr V$ which induces $(\psi, \mu)$ on $V$;
  \item the pair $(\psi, \Lambda)$ gives a Drinfeld module $\ph$ of 
    rank $2$;
  \item as the $f$-torsion of $\ph$ comes from $(\psi_f)^{-1} \Lambda/\Lambda$,
    the $f$-torsion of $\ph$ is $K_V$-rational, and thus $\ph$ has $P'(\ph)$.
\end{enumerate}

\noindent
We equip $\ph$ with a level $f$-structure as follows:  
define $\lambda(0,1) := e_{\Lambda}(\mu(1))$ and let 
$\lambda(1,0) = z$ for some element $z \in \ph[f](K_V)\bs \ph[f](V)$. 
Any $z$ gives rise to a ring homomorphism $h_z: R \lr V$ induced by 
the Weil pairing $w_f$. As before, there exists an $A_f$-automorphism
$g_z$ of $R$ with $\tilde h \circ g_z = h$, and  
$g_z$ is given by a pair $(\mf m, \s) \in {\rm Cl}(A) \times (A/fA)^*$. 
We choose an element $z$ for which $\s$ is the identity. 
As before, $z$ is unique up to a choice of $\F_q^*$. 
Clearly, the pair $(\ph, \lambda)$ is mapped to the class 
$\{ (\psi, \mu, \Lambda) \} /\F_q^*$.\\
This shows that we have a bijection between the set $$\{ \mbox{all pairs
$(\ph, \lambda)$ with $P'(\ph)$} \} / \Sigma$$ and the set
$$\{ \mbox{all triples $(\psi, \mu, \Lambda)$ with $(\psi, \mu)$ over $V$
and $(\psi_f)^{-1}\Lambda \subset K_V$} \}/\F_q^*.$$

\begin{remark}
  Similarly, the above argument shows that the injective map of Remark 
  \ref{rem_injpair} is a bijection.
\end{remark}

\par\bigskip\noindent
This bijection can be rephrased in terms of isomorphy classes of 
pairs $(\ph, \lambda)$ and triples $(\psi, \mu, \Lambda)$ as follows.  
We say that a triple $(\psi, \mu, \Lambda)$ is defined over $V$ if the
pair $(\psi, \mu)$ is defined over $V$. Two triples $(\psi, \mu,
\Lambda)$ and $(\psi', \mu', \Lambda')$ over $V$ are called {\it
isomorphic} if there exists an element $v \in V^*$ with $$(v \psi
v^{-1}, v \mu, v\Lambda) = (\psi', \mu', \Lambda').$$ Note that

\begin{equation} \label{eq_iso2}
  \{(\psi, \mu, \Lambda) \mbox{~over $V$} \}/V^* \stackrel{\rm bij}{\cong}
  \{ (\psi, \mu, \Lambda)_{K_V} \}/ K_V^*.
\end{equation}

\noindent
Moreover, if $v \in V^*$, then $$v: (\ph, \lambda) \mapsto (v
\ph v^{-1}, v \lambda)$$ and $$v: (\psi, \mu, \Lambda) \mapsto (v
\psi v^{-1}, v \mu, v \Lambda).$$ Dividing out the action of $V^*$ and
considering the bijections (\ref{eq_iso1}) and (\ref{eq_iso2}) gives
the following theorem.

\begin{theorem} \label{thm_reduct}
  Let $\Sigma = \{ \s \in \Gl_2(A/fA) \mid \det(\s) \in \F_q^* \}.$
  There is a bijection between the following two sets:
  \begin{enumerate}
    \item Isomorphy classes of pairs $(\ph, \lambda)$ over
      $K_V$ modulo $\Sigma$ where $\ph$ is a Drinfeld module
      of rank $2$ with stable reduction of rank $1$ and $\lambda$ is
      a full level $f$-structure over $K_V$.
    \item Isomorphy classes of triples $(\psi, \mu, \Lambda)$
      over $K_V$, where $\psi$ is a rank $1$ Drinfeld module over
      $K_V$, $\mu$ is a full level $f$-structure over $K_V$ and
      $\Lambda$ is an $A$-lattice of rank $1$, such that
      $(\psi_f)^{-1}\Lambda \subset K_V$.
  \end{enumerate}
\end{theorem}
\begin{proof}
  This theorem follows from the bijection that we have given above. 
\end{proof}


\section{Tate-Drinfeld modules}
\label{sec_TD}

\noindent
We follow the approach of \cite[2.2]{Boe02} and \cite{PT2} to
construct the Tate-Drinfeld module of type $\mf m$. The Tate-Drinfeld
module describes the formal neighbourhood of the cusps of the moduli scheme.
At the cusps the universal Drinfeld module with level structure degenerates
into a Drinfeld module with stable reduction. Therefore, to define
the Tate-Drinfeld module, we use the description of the stable reduction
modules as given in the previous section.\\
Let $(\psi, \mu)$ be the universal Drinfeld module of rank $1$ with
level $f$-structure over $R$. Because $f$ is invertible in $R$,
we may assume that the generator $\mu(1)$ of the $f$-torsion of $\psi$
is an invertible element in $R$. So we may and will assume that
$\mu(1) = 1$. Then by push-forward via the embeddings
$$R \lr R[\![x]\!] \lr R(\!(x)\!)$$ one has a Drinfeld module of rank $1$ with
level $f$-structure on both $R[\![x]\!]$ and $R(\!(x)\!)$. 
For an element $y = \sum_{i
\geq k} r_i x^i \in R(\!(x)\!)$ with $r_k \neq 0, k \in \Z,$ we define its
valuation $v_x(y)$ to be the $x$-valuation considered as element in
$K_R(\!(x)\!)$, and $v_x(y) = k$.

\par\bigskip\noindent
To construct a Drinfeld module of rank $2$ over $R(\!(x)\!)$, we first 
construct a lattice $\Lambda_{\mf m} \subset K_R(\!(x)\!)$,
where $\mf m \subset A$ is an ideal. This lattice turns out to depend
only on the class of $\mf m$ in the class group of $A$. As before, we can
consider $\psi_f$ as a map $$\psi_f: K_R(\!(x)\!)^{\rm sep} \lr 
K_R(\!(x)\!)^{\rm sep}.$$ The
lattice $\Lambda_{\mf m}$ will be constructed in such a way that
$(\psi_f)^{-1} \Lambda_{\mf m} \subset K_R(\!(x)\!)$. Applying Theorem
\ref{thm_reduct} to $(\psi, \mu, \Lambda_{\mf m})$ will give us the
Tate-Drinfeld module $\ph$.

\subsection{The construction of the lattice}
\noindent
Let $\mf m \subset A$ be an ideal. To prepare the construction of the lattice,
note the following:\\
1. There is a unique monic skew polynomial $P \in K_R \{ \t \}$ with
minimal degree such that $$\ker(P)(K_R^{\rm sep}) =
\psi[\mf m](K_R^{\rm sep}).$$ In fact, $P \in R \{ \t \}$
because the elements in $\ker(P)(K_R^{\rm sep})$ are integral over $R$
and $R$ is integrally closed.\\
2. Because $A_f \hookrightarrow K_R$, the extension
$K_R(\psi[\mf m])/K_R$ is Galois. Moreover, because
$\psi[\mf m](K_R^{\rm sep}) \cong A/\mf m$, there is an injective
representation $$\Gal(K_R(\psi[\mf m])/K_R) \lr (A/\mf m)^*.$$ So the Galois
group of this extension is a subgroup of $(A/\mf m)^*$.\\
3. The field $K_R(\psi[\mf m])(\!(y)\!)$ is the splitting field of the equation
$P(\frac{1}{y}) = \frac{1}{x}$ over $K_R(\!(x)\!)$. Then $$K_R(\psi[\mf
m])(\!(y)\!)/K_R(\psi[\mf m])(\!(x)\!)$$ is a Galois extension which is totally
ramified and its Galois group is isomorphic to $A/\mf m$. The
Galois action is given by $y \mapsto y + \alpha$ with $\alpha \in
\psi[\mf m](K_R(\psi[\mf m]))$.

\par\bigskip\noindent
Let $l$ be the following $A$-module homomorphism:
$$l: f^{-1} A \lr K_R(\psi[\mf m])(\!(y)\!) \quad {\rm by} \quad
f^{-1} a \mapsto \psi_a \left(\frac{1}{y} \right).$$
We use $l$ to define the lattice $\Lambda_{\mf m}$.

\begin{lemma}
  The $A$-module $l(f^{-1} \mf m)$ lies inside $R(\!(x)\!)$.
\end{lemma}
\begin{proof}
  For every $m \in \mf m$, there exists a skew polynomial $Q \in
  R \{ \t \}$ such that $\psi_m = Q \cdot P$. Note that we use here that
  one has division with remainder in the skew ring $R\{ \t \}$,
  because the leading coefficients of both $P$ and $\psi_m$ are in $R^*$.
  Consequently, $\psi_m(\frac{1}{y}) = Q(\frac{1}{x}) \in R(\!(x)\!)$.
\end{proof}

\begin{remark}
  Let $m_1, m_2$ generate $\mf m$, then there
  are elements $Q_i \in R\{ \t \}$ with $Q_i \circ P = \psi_{m_i}$. We will
  use this in the following a few times without further mentioning it.
\end{remark}

\noindent
Define the lattice $\Lambda$ as follows:
$$\Lambda := l(\mf m), \quad
{\rm then} \quad (\psi_f)^{-1}\Lambda = l(f^{-1}\mf m) + \psi[f](R)
\subset R(\!(x)\!).$$

\begin{lemma} \label{lem_latcl}
  The lattice $\Lambda$ only depends on the class of $\mf m$ in ${\rm
  Cl}(A)$.
\end{lemma}
\begin{proof}
  Suppose that $\mf m' \subset A$ is another ideal representing the
  same class as $\mf m$ in ${\rm Cl}(A)$. Then there exist elements
  $b, b' \in A$ with $b' \mf m' = b \mf m$. So we may reduce to
  case where $\mf m' = b \mf m$ for some $b \in A$. The extension
  $K_R(\psi[\mf m'])/K_R(\psi[\mf m])$ is given by adding the roots of
  the polynomial $\psi_b(X)$ to $K_R(\psi[\mf m])$.
  The extension $$K_R(\psi[\mf m'])(\!(y')\!)/K_R(\psi[\mf m'])(\!(y)\!)$$ 
  is given
  by adding the roots of the equation $\psi_b(\frac{1}{y'}) = \frac{1}{y}$.
  It is not difficult to see that $\mf m$ and $\mf m'$ give the same
  $\Lambda$.
\end{proof}

\noindent
By this lemma it makes sense to talk about the type $\mf m$ in ${\rm Cl}(A)$
of the lattice. We will denote $\mf m$ for both the ideal in $A$ as for the
ideal class in ${\rm Cl}(A)$. Moreover, this gives us some freedom in choosing
$\mf m$ to construct a lattice of type $\mf m$:

\begin{lemma} \label{lem_choice-m}
  Let $\mf a, \mf m \subset A$ be non-zero ideals, then there is an
  element $x \in K_A$, the quotient field of $A$, such that
  $x \mf m + \mf a = A$.
\end{lemma}
\begin{proof}
  By Proposition VII.5.9 in \cite{Bou}, there is an element $x \in K_A$ with
  $v_{\mf p}(x) = - v_{\mf p}(\mf m)$ for all primes $\mf p \subset A$
  dividing $\mf a$ and $v_{\mf p}(x) \geq 0$ for all other primes $\mf p$ of
  $A$. Consequently, $x \mf m \subset A$, and there is no prime ideal $\mf p$
  of $A$ dividing both $x \mf m$ and $\mf a$.
\end{proof}

\noindent
This lemma shows that we can choose a representative $\mf m$ of the class
type $[\mf m]$ of the lattice such that this representative is
relatively prime to some chosen ideal $\mf a \subset A$. This gives some
help in technical parts of some proofs later on.

\par\bigskip\noindent
We write $\Lambda_{\mf m}$ for the lattice $\Lambda$ of type $\mf m$
that we constructed above. Note that $K_R(\!(x)\!)$ is the quotient 
field of the complete discrete valuation ring $K_R[\![x]\!]$. By the 
constructions in Section \ref{sec_redth} and Theorem \ref{thm_reduct}, 
we may associate to the triple $(\psi, \mu, \Lambda_{\mf m})$ a unique Drinfeld
module $\ph$ of rank $2$ over $K_R(\!(x)\!)$ with stable reduction of rank
$1$. Moreover, the $f$-torsion of $\ph$ is $K_R(\!(x)\!)$-rational. 

\par\bigskip\noindent
In fact, $\ph$ is a Drinfeld module over $R(\!(x)\!)$, as we will now
show. 
Using Theorem \ref{thm_Dri}, the corresponding exponential map is
$$e_{\Lambda_{\mf m}}(z) := z
\prod_{\alpha \in \Lambda_{\mf m} \bs \{ 0 \}} \left( 1- \frac{z}{\alpha}
\right),$$ and $\ph$ is determined by the following diagram, which
commutes for all $a \in A$:

\begin{equation} \label{diag2}
  \begin{CD}
    0 @>>> \Lambda_{\mf m} @>>> K_R(\!(x)\!) @> e_{\Lambda_{\mf m}} >> 
    K_R(\!(x)\!) \\
    & & @VV \psi_a V @VV \psi_a V @VV \ph_a V\\
    0 @>>> \Lambda_{\mf m} @>>> K_R(\!(x)\!) @> e_{\Lambda_{\mf m}} >> 
    K_R(\!(x)\!).\\
  \end{CD}
\end{equation}

\noindent
Note that by construction $\Lambda_{\mf m} \subset R(\!(x)\!)$, each
non-zero element of $\Lambda_{\mf m}$ has negative $x$-valuation
and the leading coefficient of this non-zero element is in $R^*$.
Consequently, the map $e_{\Lambda_{\mf m}}(z) = 1 + \sum_{i\geq 1} s_i z^i$
has all its coefficients $s_i \in R[\![x]\!]$. So its inverse exists in
$R[\![x]\!][\![z]\!]$ and $$\ph_a = e_{\Lambda_{\mf m}} \circ \psi_a \circ
e_{\Lambda_{\mf m}}^{-1}.$$ Therefore, $\ph_a$ has its coefficients in 
$R[\![x]\!]$ for all $a \in A$.

\begin{lemma}
  The ring homomorphism $\ph$ is a Drinfeld module over $R(\!(x)\!)$.
\end{lemma}
\begin{proof}
  We only need to prove that the leading coefficient of $\ph_a$ is an
  element of $R(\!(x)\!)^*$. To prove this, we simply copy the computation
  of Lemma 2.10 in \cite{Boe02}.
  Note that $$\ph_a(z) = a z \prod_{\alpha \in ((\psi_a)^{-1} \Lambda_{\mf m}
  /\Lambda_{\mf m}) \bs \{ 0 \} } \left( 1 - \frac{z}{e_{\Lambda_{\mf m}}
  (\alpha)} \right).$$ So we need to show that
  $$(*) \quad \prod_{\alpha \in ((\psi_a)^{-1} \Lambda_{\mf m}
  /\Lambda_{\mf m}) \bs \{ 0 \} } e_{\Lambda_{\mf m}}(\alpha) =
  a \cdot u$$ with $u \in R(\!(x)\!)^*$\\
  Every
  $\alpha \in ((\psi_a)^{-1} \Lambda_{\mf m} /\Lambda_{\mf m}) \bs \{ 0 \}$
  which is not in $\psi[a]$ can be written uniquely as $\alpha = \alpha_1
  + \alpha_2$ where $\alpha_1$ runs through $\psi[a]$ and $\alpha_2$ runs
  through a set of representatives in 
  $\psi_a^{-1}\Lambda_{\mf m}/\Lambda_{\mf m}$ of 
  the non-zero elements of $(\psi_a)^{-1} \Lambda_{\mf m}
  /(\Lambda_{\mf m}+\psi[a])$. This set is denoted by $S_1$. The remaining
  elements $\alpha$ can be written as $\alpha = \alpha_1$ where
  $\alpha_1$ runs through $\psi[\mf m] \bs \{ 0 \}$. This set is denoted
  by $S_2$.\\
  By definition $$e_{\Lambda_{\mf m}}(\alpha) = \alpha \prod_{\beta \in
  \Lambda_{\mf m} \bs \{ 0 \}} \left( 1- \frac{\alpha}{\beta} \right).$$
  Using the rule $$\psi_a(z) = \prod_{\alpha_1 \in \psi[a]} (z - \alpha_1),$$
  we see that $$\prod_{S_1} e_{\Lambda_{\mf m}}(\alpha) =
  \prod_{\alpha_2 \neq 0} \left( \psi_a(\alpha_2) \prod_{\beta \in
  \Lambda_{\mf m} \bs \{ 0 \}} \left( \frac{\psi_a (\beta - \alpha_2)}
  {\beta^{\# A/(a)}} \right) \right).$$ This is in fact an element in
  $R(\!(x)\!)^*$, which can be seen as follows. Clearly, $0 \neq
  \psi_a(\alpha_2) \in \Lambda_{\mf m}$, so this element is in $R(\!(x)\!)^*$.
  Moreover, the element $\psi_a (\beta - \alpha_2) \in \Lambda_{\mf m}$
  cannot be $0$: if it were $0$, then any representative
  of $\alpha_2$ would lie in $\beta + \psi[a] \subset \Lambda_{\mf m} +
  \psi[a]$, i.e., the class $\alpha_2 = 0$,
  contradicting the definition of the set $S_2$. So also
  $\psi_a (\beta - \alpha_2) \in R(\!(x)\!)^*$. Finally, for almost all 
  $\beta$ we have that
  $$\frac{\psi_a (\beta - \alpha_2)}{\beta^{\# A/(a)}} \in
  R[\![x]\!]^*.$$ So the product exists, and $$\prod_{S_1}
  e_{\Lambda_{\mf m}}(\alpha) \in R(\!(x)\!)^*.$$
  On the other hand, using the rule $$h(z):= \frac{\psi_a(z)}{z} =
  \prod_{\alpha_1 \in S_2} (z - \alpha_1)$$ and recalling that $h(0) = a$,
  we see that $$\prod_{S_2} e_{\Lambda_{\mf m}}(\alpha) =
  a \cdot \prod_{\beta \in \Lambda_{\mf m} \bs \{ 0 \}}
  \left( \frac{h(\beta)}{\beta^{\# A/(a) - 1}} \right),$$ with each
  $\frac{h(\beta)}{\beta^{\# A/(a) - 1}}$ is in $R[\![x]\!]^*$.\\
  Finally, $$\prod_{\alpha \in ((\psi_a)^{-1} \Lambda_{\mf m}
  /\Lambda_{\mf m}) \bs \{ 0 \} } e_{\Lambda_{\mf m}}(\alpha) =
  \prod_{S_1} e_{\Lambda_{\mf m}}(\alpha) \cdot \prod_{S_2}
  e_{\Lambda_{\mf m}}(\alpha).$$ This finishes the proof.
\end{proof}

\par\bigskip\noindent
From the construction of the Drinfeld module $\ph$ with $\theta$ over
$R(\!(x)\!)$, coming from the triple $(\psi, \mu, \Lambda_{\mf m})$, we
can deduce at once the following list of properties:

\begin{enumerate} \label{list_tdmod}
  \item $\ph: A \lr R[\![x]\!]\{ \t \}$ is a ring homomorphism;
  \item the $f$-torsion of $\ph$ is $R(\!(x)\!)$-rational;
  \item there is an isomorphism $$A/fA \stackrel{\sim}{\lr} 
    \ph[f](R[\![x]\!])$$ given by $1 \mapsto e_{\Lambda_{\mf m}}(\mu(1))$. 
  \item $\ph~{\rm mod}~xR[\![x]\!] = \psi$, because
    $e_{\Lambda_{\mf m}}(z)~{\rm mod}~xR[\![x]\!]$ is the identity map.
\end{enumerate}

\noindent
As in Section \ref{sec_redth}, the triple $(\psi, \mu, \Lambda_{\mf m})$ 
induces on $\ph$ a level $f$-structure $\lambda$ with 
$\lambda(0,1) = e_{\Lambda_{\mf m}}(\mu(1))$ and $\lambda(1,0)$ is 
determined up to $\F_q^*$ by the Weil pairing. 

\par\bigskip\noindent
In this way, we get for every element $\mf m \in {\rm Cl}(A)$ a pair 
$(\ph^{\mf m}, \lambda^{\mf m})$. The action of $\s \in \Gl_2(A/fA)$ on 
this pair is given by $$\s: (\ph^{\mf m}, \lambda^{\mf m}) \lr 
(\ph^{\mf m}, \lambda^{\mf m} \circ \s).$$
The action of the class group of $A$ is described in the following lemma.
Let $\mf n \in {\rm Cl}(A)$ and let $g_{\mf n}$ denote the $A_f$-linear
automorphism of $R$ which describes the action of $\mf n$ on $R$.
Let $\Lambda_{\mf n^{-1} \mf m}$ denote the lattice of type 
$\mf n^{-1}\mf m$.

\begin{lemma} \label{lem_actClA}
  The element $\mf n$ maps the triple $(\psi, \mu, \Lambda_{\mf m})$ 
  to the triple $$(g_{\mf n}(\psi), g_{\mf n}(\mu), 
  g_{\mf n}(\Lambda_{\mf n^{-1}\mf m})).$$ 
\end{lemma}
\begin{proof}
  We choose representatives $\mf m, \mf n \subset A$ of the classes 
  of $\mf m$ and $\mf n$ in ${\rm Cl}(A)$ such that 
  $\mf n^{-1} \mf m \subset A$, and $\mf n$ and $\mf m$ are relatively
  prime to $f$.\\ 
  Write $\ph = \ph^{\mf m}$.
  The action of $\mf n$ on $\ph$ is given by a unique monic skew
  polynomial $Q$ with minimal degree such that
  $$\ker(Q)(K_R(\!(x)\!)^{\rm sep}) = \ph[\mf n](K_R(\!(x)\!)^{\rm
  sep}).$$ Let $\ph'$ be the image under $Q$: $$\ph
  \stackrel{Q}{\lr} \ph'.$$ Writing $\Lambda_{\mf m} = l(\mf m)$
  as before, it is not difficult to see that
  $$\ker(Q\circ e_{\Lambda_{\mf m}})(K_R(\!(x)\!)^{\rm sep}) =
  l(\mf n^{-1} \mf m) + \psi[\mf n](K_R(\!(x)\!)^{\rm sep}).$$
  The action of $\mf n$ on $R$, denoted by $g_{\mf n}$, corresponds to 
  a skew polynomial $P \in R\{ \t \}$ with $\ker(P)(K_R^{\rm sep}) = 
  \psi[\mf n](K_R^{\rm sep})$. Let $\Lambda'$ be the 
  lattice given by $$\Lambda':= P\circ 
  \ker(Q\circ e_{\Lambda_{\mf m}})(K_R(\!(x)\!)^{\rm sep}) =
  P \circ l(\mf n^{-1} \mf m) $$ $$= \{ g_{\mf n}(\psi_{fa})
  P(\frac{1}{y}) \mid a \in \mf n^{-1} \mf m \}.$$ Then $\Lambda' = 
  g_{\mf n}(\Lambda_{\mf n^{-1}\mf m})$, and it is not difficult to see 
  that $\ph'$ corresponds to the pair $(g_{\mf n}(\psi), 
  g_{\mf n}(\Lambda_{\mf n^{-1}\mf m}))$. 
\end{proof}

\par\bigskip\noindent
As a corollary to this lemma, we see that $\mf n$ maps $\ph^{\mf m}$ to 
$g_{\mf n}(\ph^{\mf n^{-1} \mf m})$. As the definition of $\lambda^{\mf m}$
and $\lambda^{\mf n^{-1} \mf m}$ depend on a choice of $\F_q^*$, we cannot
say that $\mf n$ maps $(\ph^{\mf m}, \lambda^{\mf m})$ to 
$(\ph^{\mf n^{-1} \mf m}, \lambda^{\mf n^{-1} \mf m})$. 
Bearing this in mind, we propose for every type $\mf m$ the following 
definition of the pairs $(\tdm, \ldm)$. For $\mf m = A$ we define 
$(\tdm, \ldm) := (\ph^{\mf m}, \lambda^{\mf m}).$ 
This pair is unique up to $\F_q^*$. For the rest of this paper we keep 
this choice is fixed. In general, $\mf m^{-1}$ gives rise 
to an automorphism $g_{\mf m^{-1}}$ of $R$. We extend this to an automorphism 
of $R[\![x]\!]$ by letting it act trivially on $x$. We define 
$(\tdm, \ldm)$ to be the image of $(\td_{A}, \ld_{A})$ under $\mf m^{-1}$.

\begin{definition} \label{def_TDlev}
  The pair $(\tdm, \ldm)$ is called the {\em standard Tate-Drinfeld 
  module of rank $2$ and type $\mf m$ with level $f$-structure}.
  A pair $(\tdm, \lambda)$, where $\lambda$ is any level $f$-structure
  over $R(\!(x)\!)$ is called a {\it Tate-Drinfeld module with level 
  $f$-structure.}
\end{definition}

\noindent
The pairs $(\xi \tdm \xi^{-1}, \xi \lambda)$ with $\xi \in R[\![x]\!]^*$
all belong to the same isomorphy class $(\tdm, \lambda)$. 
Every isomorphy class of Tate-Drinfeld modules with level $f$-structure
comes from a unique morphism $$\Spec(R(\!(x)\!)) \lr M^2(f).$$
Because we always assume that $\mu(1) = 1$, the pair $(\tdm, \ldm)$ is fixed 
in its isomorphy class. For every $\lambda$ there is a unique $\s \in 
\Gl_2(A/fA)$ with $\lambda = \ldm \circ \s$. Therefore, the pair 
$(\tdm, \lambda)$ is fixed in its isomorphy class.


\section{The universal Tate-Drinfeld module}
\label{sec_UP}

\noindent
In this section we introduce the universal Tate-Drinfeld module. 
This is a scheme $\Cal Z$ consisting of the coproduct of a number
of copies of $\Spec(R[\![x]\!])$ and which is equipped with a 
Tate-Drinfeld structure $(\td, \ld)$. Let $N \subset
\Gl_2(A/fA)$ be the subgroup $$N = \left( \begin{array}{cc} \F_q^* & A/fA \\
0 & (A/fA)^* \end{array} \right).$$ Choose representatives
$\s_1, \ldots, \s_n$ of the cosets of $N \s_i$ in 
$\Gl_2(A/fA)$ and let $\s_1$ be the identity. As $\det(N) = (A/fA)^*$, 
we may assume that the representatives $\s_i$ are elements of $\Sl_2(A/fA)$.
Set $$\Cal Z = \Spec \left( \underset{(\mf m, \s_i)}{\oplus}
R[\![x]\!]_{(\mf m, \s_i)} \right),$$
The pair $(\td, \ld)$ is the Drinfeld module $\td$ of rank $2$ with 
level $f$-structure $\ld$ over $\Spec \left( \underset{(\mf m, \s_i)}{\oplus}
R(\!(x)\!)_{(\mf m, \s_i)} \right)$ such that the
restriction of $(\td, \ld)$ to $R(\!(x)\!)_{(\mf m, \s_i)}$ is equal to 
$(\tdm, \ldm \circ \s_i)$ for every pair $(\mf m, \s_i)$.  

\begin{definition}
  The {\it universal Tate-Drinfeld module of rank $2$ with level
  $f$-structure} is the scheme $\Cal Z$ together with the pair $(\td, \ld)$.
\end{definition}

\begin{remark} \label{rem_Act}
  As before, a pair $(\mf m, \s_i) \in {\rm Cl}(A)
  \times \Gl_2(A/fA)$ determines a unique action on $(\td_{A}, \ld_{A})$. 
  And by definition we have $$(\tdm, \ldm \circ \s_i) = (\mf m^{-1}, \s_i) 
  (\td_{A}, \ld_{A}).$$
  The action of ${\rm Cl}(A) \times \Gl_2(A/fA)$ on the universal
  Tate-Drinfeld module 
  is given by this action of $(\mf m, \s_i)$ and the fact that $N$ acts 
  trivially. 
\end{remark}

\par\bigskip\noindent
Clearly, we would like to have a certain universal property for this universal
Tate-Drinfeld module. The weak versions of this universal property that
we need can be found in Theorems \ref{thm_univTD} and
\ref{thm_STD_up}. In Proposition \ref{prop_STD_up} the main work is done
for Theorem \ref{thm_STD_up}. This proposition explains the
subgroup $N$.

\begin{prop} \label{prop_STD_up}
  Let $\s \in \Gl_2(A/fA)$. There exists an $A_f$-linear ring homomorphism
  $h_{\s}: R[\![x]\!] \lr R[\![x]\!]$ such that $$h_{\s}(\tdm, \ldm) \cong
  (\tdm, \ldm \circ \s)$$ if and only if $\s \in N$. 
  If $\s \in N$, then $h_{\s}$ is given by $\det(\s): R \lr R$ and 
  $x \mapsto \delta x$ for some $\delta \in R[\![x]\!]^*$.
  Moreover, for any $\s \in N$ the image of $(\tdm, \ldm \circ \s_i)$
  under $h_{\s}$ is isomorphic to $(\tdm, \ldm \circ \s \circ \s_i)$. 
\end{prop}

\begin{proof}
  We write $\s = (\s_{i,j})$. First suppose that $h_{\s}$ exists. Because
  $h_{\s}$ is determined by $\det(\s)$ and $h_{\s}(x)$, it must respect
  the ordering of the $x$-valuation. This implies that
  $h_{\s} (\lambda(1,0))$ must
  have minimal negative valuation and $h_{\s} (\lambda(0,1))$ must
  have valuation $0$. Hence the only $\s$'s whose action may come from an
  $A_f$-linear ring homomorphism $h_{\s}$ are $\s \in N$. This shows the
  'only if'.\\
  We prove the 'if'-part in two steps. Let $\s = (\s_{i,j}) \in N$.\\
  1. First suppose that $\s = \left( \begin{array}{cc} 1 & 0 \\ 0 & \alpha
  \end{array} \right)$ with $\alpha \in (A/fA)^*$. The action of $\s$
  induces an action $$\alpha = \det(\s): M^1(f) \lr M^1(f),$$ i.e., there
  is an $A_f$-linear ring homomorphism $h_{\alpha}: R \lr R$ such that
  $$(h_{\alpha}(\psi), h_{\alpha}(\mu)) = (\xi \psi \xi^{-1}, \xi \mu
  \alpha) \cong (\psi, \mu \alpha),$$ for some $\xi \in R^*$.
  We use the notations of Section \ref{sec_TD}. Recall that
  the element $\frac{1}{y}$ was used to define the lattice
  $\Lambda^{\rm td}_{\mf m}$ and that $\frac{1}{x} = P(\frac{1}{y})$. 
  $P \in R\{ \t \}$ is the skew polynomial of minimal degree with $$\ker(P)
  (K_R^{\rm sep}) = \psi[\mf m](K_R^{\rm sep}).$$ The map $\frac{1}{y}
  \mapsto \xi \frac{1}{y}$ induces, by applying $P$, $$\frac{1}{x} \mapsto
  h_{\alpha}(P)(\xi \frac{1}{y}) = \delta^{-1} \frac{1}{x} \quad \mbox{for
  some $\delta^{-1} \in R[\![x]\!]^*$.}$$ To see this, note that $h_{\alpha}(P)
  = \zeta P \xi^{-1}$ for some $\zeta \in R[\![x]\!]^*$. (In fact, $\zeta$
  is determined by the fact that $h_{\alpha}(P)$ is monic.)\\
  The map $h_{\alpha}$ is extended to a ring homomorphism 
  $$h_{\s}: R[\![x]\!] \lr R[\![x]\!], \quad x \mapsto \delta x,
  \quad h_{\s}(r) = h_{\alpha}(r) \quad \forall r \in R.$$
  An easy computation shows that $$h_{\s}(\Lambda^{\rm td}_{\mf m}) = \xi
  \Lambda^{\rm td}_{\mf m},\quad \mbox{and thus} \quad h_{\s}(\tdm) = \xi
  \tdm \xi^{-1}.$$ Using that there is an element $m \in \mf m$ with
  $\lambda(1,0) = e_{\Lambda^{\rm td}_{\mf m}}(\psi_m(\frac{1}{y}))$, we see 
  that 
  $$h_{\s}
  \left( \begin{array}{c} \lambda(1,0) \\ \lambda(0,1) \end{array} \right) =
  \left( \begin{array}{c} e_{\xi \Lambda^{\rm td}_{\mf m}}(\xi
  \psi_m(\frac{1}{y})) \\ e_{\xi \Lambda^{\rm td}_{\mf m}}(\xi \mu(1) \alpha)
  \end{array} \right) =
  \left( \begin{array}{c} \xi \lambda(1,0) \\ \xi \lambda(0,1) \alpha
  \end{array} \right).$$
  And so indeed, $$h_{\s}(\tdm, \lambda) = (\xi \tdm \xi^{-1},
  \xi \lambda \circ \s) \cong (\tdm, \lambda \circ \s).$$
  2. Having dealt with the first case, we may assume that $\s \in N$ and
  $\det(\s) \in \F_q^*$. As a consequence, the map $h_{\s}$ that we are looking
  for is $R$-linear and in particular, if $(\psi, \mu, \Lambda)$ is the
  triple associated to $(\tdm, \lambda)$, then $h_{\s}(\psi) = \psi$
  and $h_{\s}(\mu) = \mu$.\\ As always, the action of $\F_q^* \cdot \left(
  \begin{array}{cc} 1 & 0 \\ 0 & 1 \end{array} \right)$ is trivial, so we
  may assume that $\s_{2,2} = 1$.
  We first give a proof in the simple case $\Lambda \cong A$ and then prove
  it for general $\Lambda$.\\
  i. Suppose $\Lambda \cong A$. Note that in that case, the elements
  $e_{\Lambda}(\frac{1}{x})$ and $e_{\Lambda}(\mu(1))$ generate the 
  $f$-torsion
  of $\td$. In fact, there is a basis transformation $$\alpha = \left(
  \begin{array}{cc} \alpha_{1,1} & \alpha_{1,2} \\ 0 & 1 \end{array} \right)
  \quad \mbox{with} \quad \alpha
  \left( \begin{array}{c} \lambda(1,0) \\ \lambda(0,1) \end{array} \right) =
  \left( \begin{array}{c} e_{\Lambda}(\frac{1}{x}) \\ e_{\Lambda}(\mu(1))
  \end{array} \right)$$ and with $\alpha_{1,1} \in (A/fA)^*$, 
  $\alpha_{1,2} \in A/fA$.\\
  Let $\rho = \alpha \s \alpha^{-1}$, then $\rho = \left(
  \begin{array}{cc} \rho_{1,1} & \rho_{1,2} \\ 0 & 1 \end{array} \right)$
  with $\rho_{1,1} \in \F_q^*$, $\rho_{1,2} \in A/fA$.
  We set $$\delta^{-1} = \rho_{1,1} + \rho_{1,2} (\mu(1)) \cdot x \in
  R[\![x]\!]^*,$$ and we define $h_{\s}$ to be the $R$-linear map given by
  $$h_{\s}: x \mapsto \delta \cdot x.$$ Because $$\psi_f(\delta^{-1}
  \frac{1}{x}) = \rho_{1,1} \psi_f(\frac{1}{x}),$$ it follows that
  $h_{\s}(\Lambda) = \Lambda$, and thus $h_{\s}$ commutes with $e_{\Lambda}$.
  Moreover $$h_{\s} (\psi, \mu, \Lambda) = (\psi, \mu, \Lambda).$$ Therefore
  $$h_{\s}(\tdm) = \tdm.$$
  So $h_{\s}$ commutes with the $A$-action and with $e_{\Lambda}$. To see
  what happens on the level structure $\lambda$ is now an easy computation.
  $$h_{\s}
  \left( \begin{array}{c} \lambda(1,0) \\ \lambda(0,1) \end{array} \right) =
  h_{\s} \alpha^{-1}
  \left( \begin{array}{c} e_{\Lambda}(\frac{1}{x}) \\ e_{\Lambda}(\mu(1))
  \end{array} \right) =$$
  $$\alpha^{-1} \rho \left( \begin{array}{c} e_{\Lambda}(\frac{1}{x}) \\
  e_{\Lambda}(\mu(1)) \end{array} \right) =
  \s \left( \begin{array}{c} \lambda(1,0) \\ \lambda(0,1) \end{array}
  \right).$$
  We conclude that $$h_{\s}(\ph^{\rm td}, \lambda^{\rm td}) =
  (\ph^{\rm td}, \lambda^{\rm td} \circ \s).$$
  ii. In general, let $\Lambda \cong \mf m$ and recall the construction of
  $\Lambda$ in Section \ref{sec_TD}. Using the same notations as in Section
  \ref{sec_TD}, let $m \in \mf m$ be an element such that the image of
  $\frac{1}{z} := \psi_m(\frac{1}{y})$ and $\mu(1)$ under $e_{\Lambda}$
  generate the $f$-torsion of $\tdm$. We let $\alpha$ be as in i,
  i.e., $$\alpha \left( \begin{array}{c} \lambda(1,0) \\ \lambda(0,1)
  \end{array} \right) = \left( \begin{array}{c} e_{\Lambda}(\frac{1}{z}) \\
  e_{\Lambda}(\mu(1)) \end{array} \right)$$ and $\rho = \alpha \s \alpha^{-1}$.
  As in the previous case, we look for an $h_{\s}$ such that
  $$\frac{1}{z} \mapsto (\rho_{1,1} + \rho_{1,2}
  \mu(1) z) \cdot \frac{1}{z}.$$ This can be done as follows.
  We start by assuming that $m \in (A/fA)^*$ - we will show below that we
  may assume this in general. Let $b \in A$ such that $b \equiv m^{-1}
  \rho_{1,2} \in A/fA$. Let $$\frac{1}{y} \mapsto
  \rho_{1, 1} \frac{1}{y} + \psi_{b} \mu(1).$$ Applying $P$  gives our
  candidate for $h_{\s}$: $$\frac{1}{x} \mapsto
  \delta^{-1} \frac{1}{x} \quad ~{\rm with}~
  \delta^{-1} = \rho_{1, 1} + P(\psi_{b} (\mu(1))) \cdot x.$$ 
  By construction
  $$h_{\s}: \frac{1}{z} \mapsto \rho_{1,1} \frac{1}{z} +
  \psi_m(\psi_b (\mu(1))) = \rho_{1,1} \frac{1}{z} + \rho_{1,2} \mu(1).$$
  Note that all elements of $\Lambda$ have the form $\psi_{f \tilde m}
  (\frac{1}{y})$, with $\tilde m \in \mf m$. And clearly,
  $$h_{\s} (\psi_{f \tilde m}(\frac{1}{y})) = \rho_{1,1} \psi_{f \tilde m}
  (\frac{1}{y}).$$ Because $\rho_{1,1} \in \F_q^*$, we see that
  $h_{\s}(\Lambda) = \Lambda$.\\
  We can conclude the proof in the same way as in the case of $\Lambda
  \cong A$.\\
  Finally, it remains to be shown that $m$, which
  we used to define $\frac{1}{z} = \psi_m(\frac{1}{y})$, is an element of
  $(A/fA)^*$. By Lemma \ref{lem_choice-m} we may assume that $\mf m$ and
  $(f)$ are relatively prime. Furthermore, because $\frac{1}{z}$ generates
  one direct summand of the $f$-torsion, one has $\psi_b(\frac{1}{z})
  \in \Lambda$ if and only if $b \in fA$, i.e., $b m \in f \mf m$ if and only
  if $b \in fA$. Consequently, $(f) + (m) = A$.\\
  The `moreover'-part of the proposition is obvious. 
\end{proof}

\noindent
Using Proposition \ref{prop_STD_up}, we can immediately prove one weak form
of the universal property of the universal Tate-Drinfeld module.

\begin{theorem} \label{thm_STD_up}
  For every Tate-Drinfeld module $(\tdm, \lambda)$ there is a unique
  ring homomorphism $$h: \underset{(\mf m', \s_i)}{\oplus} R[\![x]\!]_{(\mf m',
  \s_i)} \lr R[\![x]\!]$$
  such that $$h(\td, \ld) \cong (\tdm, \lambda).$$
\end{theorem}
\begin{proof}
  Let $\s \in \Gl_2(A/fA)$ such that $\ldm \circ \s = \lambda$
  and let $\s_k \in \{ \s_1, \ldots, \s_n \}$ such that $\s \in \s_k N$.
  The map $h$ is defined as follows: $h$ is the zero-map on 
  $R[\![x]\!]_{(\mf m', \s_i)}$ if $\mf m' \neq \mf m$ in the class group 
  of $A$ or if
  $\s_k \neq \s_i$. On $R[\![x]\!]_{(\mf m, \s_k)}$ it is the map defined
  in Proposition \ref{prop_STD_up}.
  To show uniqueness, note that any $A_f$-linear ring homomorphism $$h:
  R[\![x]\!]_{(\mf m', \s_i)} \lr R[\![x]\!]$$ induces on $R[\![x]\!]$ a 
  Tate-Drinfeld
  module whose corresponding lattice has type $\mf m'$. Hence,
  any such ring homomorphism keeps the type of the Tate-Drinfeld module fixed.
  Moreover, if there is an $A_f$-morphism
  $$R[\![x]\!]_{(\mf m, \s_i)} \lr R[\![x]\!]_{(\mf m, \s_j)}$$ which 
  induces the Tate-Drinfeld structure, then there is a morphism
  $$R[\![x]\!]_{(\mf m, \s_1)} \lr R[\![x]\!]_{(\mf m, \s_j \s_i^{-1})},$$ 
  and thus by Proposition \ref{prop_STD_up} we see that
  $\s_j \s_i^{-1} \in N$. So $\s_i = \s_j$.
\end{proof}

\subsection{The weak version of the universal property of $\Cal Z$}
\noindent
Let, as in Section \ref{sec_redth}, $V$ be a complete discrete valuation
$A_f$-ring, let $\pi$ be a generator of its maximal ideal, and let $K_V$ 
be its field of
fractions. Let $(\ph, \lambda)$ be a Drinfeld module $\ph$ of rank $2$ over
$K_V$ with level $f$-structure $\lambda$ such that $\ph$ has
stable reduction of rank $1$ at $\pi$. In this subsection we
discuss the other weak version of the universal property of the
universal Tate-Drinfeld module, which we
need in the next section. We prove that there exists a unique ring
homomorphism
$$h_{\ph, \lambda}: \underset{(\mf m, \s_i)}{\oplus} R[\![x]\!] \lr V$$
such that $$h_{\ph, \lambda}(\td, \ld) \cong (\ph, \lambda).$$

\par\smallskip\noindent
In Theorem \ref{thm_reduct}
we showed that each triple $(\ph, \lambda, z)$ corresponds modulo
$\Sigma$ to a unique triple $(\psi, \mu, \Lambda)$. Let $\mf m$
be the type of $\Lambda$. In trying to avoid confusion, we will write
$(\psi^{\rm un}, \mu^{\rm un}, \Lambda_{\mf m}^{\rm td})$ for the triple
used in defining the Tate-Drinfeld module of type $\mf m$ over $R(\!(x)\!)$,
where indeed $(\psi^{\rm un}, \mu^{\rm un})$ is the universal Drinfeld
module with level $f$-structure of rank $1$ over $R$.\\
The pair $(\psi, \mu)$ over $V$ comes from an $A_f$-linear ring 
homomorphism, which we called $\tilde h$: $$\tilde h: R \lr V.$$ We have 
$\tilde h(\psi^{\rm un}, \mu^{\rm un}) =(\psi, \mu)$. 

\par\bigskip\noindent
We show that there exists an extension of $\tilde h$ to $R[\![x]\!]$
$$h_{\psi, \mu}: R[\![x]\!] \lr V,$$ where $R[\![x]\!]$ comes equipped 
with $\tdm$, such that $h_{\psi, \mu}(\tdm) = \ph$; cf. Proposition
\ref{prop_stab-from-td}. The main point is showing that there exists a 
ring homomorphism $h_{\psi, \mu}$ such that
$$h_{\psi, \mu}((\psi^{\rm un}_f)^{-1}
\Lambda^{\rm td}_{\mf m}) = (\psi_f)^{-1} \Lambda.$$

\begin{lemma} \label{lem_gen-elm-lat}
  Let $m_1$ and $m_2$ generate the ideal $\mf m$. Then there
  exists an element $\zeta \in \ov K_V$ such that the projective part of
  $(\psi_f)^{-1} \Lambda$ is generated as $A$-module by the elements
  $\psi_{m_1}(\zeta)$ and $\psi_{m_1}(\zeta)$.
\end{lemma}
\begin{proof}
  Let $M = \ov K_V$ be the algebraic closure of $K_V$, and let $M_{\rm
  tor}$ be the set of $A$-torsion points in $M$; the $A$-action
  is given by $\psi$.\\
  1. The $A$-module $M_{\rm tor}$ is divisible, i.e., for all $a \in A
  \backslash \{ 0 \}$ the map $$\psi_a: M_{\rm tor} \lr M_{\rm tor}$$
  is surjective. Namely, if $x \in M_{\rm tor}$, then the equation
  $\psi_a(z) = x$ has solutions in $M_{\rm tor}$. Consequently, 
  $M_{\rm tor}$ is an
  injective module. Cf. Theorem 7.1 in \cite{HS97}], where it is shown that 
  divisibility is the same as injectivity for modules over a PID; it
  is not difficult to extend this to a theorem over Dedekind domains.\\
  2. $M/M_{\rm tor}$ has a natural $K_A$-module structure. 
  The following sequence of $A$-modules is
  exact: $$0 \lr M_{\rm tor} \lr M \lr M/M_{\rm tor} \lr 0.$$ Note
  that $\Lambda$ is torsion-free, hence $\Lambda \oplus
  M_{\rm tor} \hookrightarrow M$. Consider the projection map $$s:
  \Lambda \oplus M_{\rm tor} \lr M_{\rm tor}, \quad {\rm by}~\alpha
  \oplus m \mapsto m.$$ Because $M_{\rm tor}$ is an injective module, it
  follows that $s$ extends to a map $$s: M \lr M_{\rm tor}.$$ So the
  exact sequence splits according to $s$ and $M \cong M_{\rm tor}
  \oplus M/M_{\rm tor}$.\\
  3. According to $2.$, we may write $(\psi_f)^{-1} \Lambda =
  N_1 \oplus N_2$ with $N_1 = \psi[f](V)$ being the torsion part
  and $N_2 \cong \mf m$ being the projective part of $(\psi_f)^{-1}\Lambda$.
  Let $e_1, e_2$ be generators of $N_2$ such that $\psi_{m_2} e_1 =
  \psi_{m_1} e_2$. Then $$(\psi_{m_1})^{-1} e_1 \equiv (\psi_{m_2})^{-1}
  e_2 \mod M_{\rm tor}.$$ Let $\zeta_i \in M$ be the unique
  element with $\zeta_i \mapsto (\psi_{m_i})^{-1} e_i \in M/M_{\rm
  tor}$ and $s(\zeta_i) = 0$. Then $\zeta_1 - \zeta_2 \mapsto 0 \in
  M/M_{\rm tor}$ and $s(\zeta_1 - \zeta_2) = 0$. Consequently,
  $\zeta:= \zeta_1 = \zeta_2$ is the element we are looking for.
\end{proof}

\begin{prop} \label{prop_stab-from-td}
  Every rank $2$ Drinfeld module $\ph$ over $V$ with $K_V$-rational
  $f$-torsion and with stable reduction of rank $1$ and type $\mf m$ is
  induced by $\tdm$ via the ring homomorphism 
  $$h_{\psi, \mu}: R[\![x]\!] \lr V,$$
  i.e., $h_{\psi, \mu}(\tdm) \cong_V \ph$. 
  Extending $h_{\psi, \mu}$ to $$R(\!(x)\!) \lr K_V$$ maps the $f$-torsion
  of $\tdm$ isomorphically to the $f$-torsion of $\ph$.
\end{prop}
\begin{proof}
\noindent
  In Section \ref{sec_TD} we introduced skew polynomials
  $P, Q_1, Q_2 \in R\{ \t \}$ with $Q_i \circ P = \psi^{\rm un}_{m_i}$. 
  In particular, we see that $\tilde h(P)$ divides $P':= {\rm gcd}( 
  \psi_{m_1}, \psi_{m_2}) \in K_V\{ \t \}$. We may assume 
  that $\mf m$ is not contained in the kernel of $A_f \lr V$. Therefore, 
  $\deg_{\t} P' = \deg_{\t} h_{\psi, \mu}(P)$. Consequently, there exist
  elements $\beta_i \in K_V \{ \t \}$ with
  $$\tilde h(P) = \beta_1 \psi_{m_1} + \beta_2 \psi_{m_2}.$$
  Let $\zeta$ be the element from Lemma \ref{lem_gen-elm-lat}, and 
  define $\frac{1}{z} := \tilde h(P)(\zeta).$
  As $\psi_{m_1}(\zeta)$ and $\psi_{m_1}(\zeta)$ generate the projective
  part of $(\psi_f)^{-1} \Lambda \subset K_V$, we see that $\frac{1}{z} \in 
  K_V$ and $z \in V$.\\
  We extend $\tilde h$ to 
  $$h_{\psi, \mu}: R[\![x]\!] \lr V \quad \mbox{by} \quad x \mapsto z.$$
  One can easily verify that $$h_{\psi, \mu}((\psi^{\rm un}_f)^{-1}
  \Lambda^{\rm td}_{\mf m}) = (\psi_f)^{-1} \Lambda.$$
\end{proof}

\begin{remark}
  The proof of Proposition 2.5 in \cite{PT2} is not entirely complete.
  One way of completing it, is by adding the construction of the lattice as 
  is done here. This makes sure that the morphism $h$ in Proposition 2.5 in 
  \cite{PT2} indeed exists. 
\end{remark}

\begin{theorem} \label{thm_univTD}
  Every pair $(\ph, \lambda)$ consisting of a rank $2$ Drinfeld module 
  $\ph$ over $V$ with stable reduction of rank $1$ with level $f$-structure
  $\lambda$ is induced by $(\td, \ld)$ via the unique ring homomorphism
  $h_{\ph, \lambda}$.
\end{theorem}
\begin{proof}
  Let $(\psi, \mu, \Lambda)$ be the triple associated to $(\ph, \lambda)$, 
  and let $\mf m$ be the type of $\Lambda$. Let $h_{\psi, \mu}: R[\![x]\!] 
  \lr V$ be the morphism such that $h_{\psi, \mu}(\tdm) = \ph$. 
  Let $\s \in \Gl_2(A/fA)$ such that $h_{\psi, \mu}(\ldm \circ \s) = \lambda$.
  The element $\s$ lies in a unique class $N \s_i$. Write $\s = \t \s_i$ 
  with $\t \in N$. Let $h_{\t}$ be the ring homomorphism which is defined in 
  Proposition \ref{prop_STD_up}.
  Define $$h_{\ph, \lambda}: \oplus R[\![x]\!]_{(\mf m', \s_j)} \lr V$$ as
  follows: $h_{\ph, \lambda}$ equals $h_{\psi, \mu} \circ h_{\t}$ on
  $R[\![x]\!]_{(\mf m, \s_i)}$ and is zero on the other copies of 
  $R[\![x]\!]$.\\
  The uniqueness follows from the construction and from Theorem 
  \ref{thm_STD_up}.
\end{proof}


\section{The compactification of $M^2(f)$} \label{sec_compact}

\noindent
In this section we describe a compactification $\ov M^2(f)$ of
$M^2(f)$, which is analogous to the compactification of the classical
modular curves given by Katz and Mazur in Chapter 8 of \cite{KM85}. We
define the scheme of cusps, which we call {\it Cusps}. This is a
closed subscheme of $\ov M^2(f)$. Moreover, we consider the formal
scheme $\widehat{\it Cusps}$, which is the completion of $\ov M^2(f)$
along {\it Cusps}. In the following section we will use the universal
Tate-Drinfeld module as defined in the previous sections to describe
the scheme of cusps.

\subsection{The morphism $j_a$}
\noindent
Let $a \in A \bs \F_q$. 
Let $(\ph, \lambda)$ be the universal Drinfeld module of rank $2$ with
level $f$-structure over $M^2(f)$. Let $B$ be the ring with $\Spec(B)
= M^2(f)$ and write $$\ph_a = \sum_{i=0}^{2\deg(a)} b_i \t^i \quad
\mbox{with $b_i \in B$ for all $i$ and $b_{2\deg(a)} \in B^*$}.$$
Let $j_a: M^2(f) \lr \A^1_{A_f}$ be the morphism given by 
$$j_a^{\#}: A_f[j] \lr B, \quad j \mapsto b_{\deg(a)}^{q^{\deg(a)}
+ 1}/b_{2\deg(a)};$$ cf. \cite[4.2]{Leh00}. Clearly, $j_a$ factors 
over $M^2(1)$.

\begin{lemma} \label{lem_ja}
  The morphism $j_a$ is finite and flat.
\end{lemma}
\begin{proof}
  (This is the proof of Proposition 4.2.3 in \cite{Leh00}.) The morphism
  $j_a$ is of finite type, hence we may use the valuative criterion to
  prove properness. Suppose that $V$ is a discrete valuation ring, let
  $K_V$ be its quotient field, and suppose that there are morphisms given
  that make the following diagram commutative.
  $$\begin{CD} \Spec(K_V)
  @>>> M^2(f)\\ @VVV @VV j_a V\\ \Spec(V) @>>> \A^1_{A_f}.
  \end{CD}$$
  Note that the upper horizontal map gives rise to a (unique) map
  $$\Spec(V) \lr M^2(f)$$ if and only if the pull-back of the pair $(\ph,
  \lambda)$ via the upper horizontal map has good reduction at the
  maximal ideal of $V$. So $j_a$ is proper if and only if $j_a(x)
  \in V$ implies that the pull-back $(\ph', \lambda')$ over $K_V$ has
  good reduction. If this
  pull-back does not have good reduction, it has stable reduction of
  rank $1$. Suppose that $(\ph, \lambda)$ has stable reduction of
  rank $1$, then there exists an element $s\in K_V^*$ such that $s
  \ph_a s^{-1}$ has all coefficients in $V$, the
  $\deg(a)$'th coefficient has valuation $0$, and the $2\deg(a)$'th
  coefficient has strictly positive valuation. This means that the image of
  $j_a(x)$ is not in $V$. We conclude that $j_a$ is proper.\\
  Because each connected component of $M^2(f)$ is an affine variety
  over $\F_q$, it follows by \cite[ex. II.4.6]{Har77} that $j_a$
  restricted to such a connected component is finite. And thus $j_a$
  is finite.\\ The finite ring homomorphism $j_a^{\#}$ is
  injective. Let $\mf P \subset B$ be a prime ideal lying above $\mf p
  \subset A_f[j]$. Then both local rings $B_{\mf P}$ and $A_f[j]_{\mf
  p}$ are regular and of equal dimension. By the finiteness of $j_a$
  it follows that $B_{\mf P}$ is a free $A_f[j]_{\mf p}$-module, cf. 
  Corollary IV.22 in \cite{Ser65}. Hence, $j_a$ is flat.
\end{proof}

\subsection{The compactification}
\noindent
The ring $B$ is a finite $A_f[j]$-algebra via $j_a^{\#}$. 
Let $C$ denote the normalization of $A_f[\frac{1}{j}]$ inside the 
quotient ring of $B$. Then $C$ is finite over $A_f[\frac{1}{j}]$;
cf. Corollary 13.13 in \cite{Eis95}.\\
The compactification $\ov M^2(f)$ of $M^2(f)$ is defined as the 
scheme obtained by glueing $\Spec(B)$ and $\Spec(C)$ along their
intersection. We obtain a finite morphism 
$$\ov j_a: \ov M^2(f) \lr \P^1_{A_f}.$$ The following diagram is cartesian
$$\begin{CD}
  M^2(f) @>>> \ov M^2(f)\\
  @V j_a VV @VV \ov j_a V\\
  \A^1_{A_f} @>>> \P^1_{A_f}.
\end{CD}.$$

\begin{remark}
  By Theorem \ref{thm_CompToTD} it follows that $\ov j_a$ is flat
  in the points of the boundary of $\ov M^2(f)$; therefore, $\ov j_a$ is 
  flat. 
\end{remark}

\begin{lemma} \label{lem_ind_comp}
  The scheme $\ov M^2(f)$ is independent of the chosen element $a$.
\end{lemma}
\begin{proof}
  Let $B'$ be a connected component of $B$.  
  Let $a_1, a_2 \in A \bs \F_q$ and consider the maps 
  $j_{a_i}: A_f[j_i] \lr B'$. Let $C_i$ be the integral closure 
  of $A_f[\frac{1}{j_i}]$ inside $K_{B'}$. Let $X_i$ be the scheme 
  obtained by glueing $\Spec(C_i)$ and $\Spec(B')$ along their 
  intersection.\\ 
  Let $\mf p \in X_1 \bs \Spec(B')$ be any prime of height one of 
  $C_1$, then $\frac{1}{j_1} \in \mf p$. Let $v_{\mf p}$ be the valuation 
  of $K_{B'}$ given by $\mf p$. If $v_{\mf p}(\frac{1}{j_2}) \leq 0$, then 
  $\mf p$ would correspond to a valuation of $K_A[j_2]$ and therefore to
  a valuation of $B'$; cf. Section VII.9 in \cite{Bou}. 
  Consequently, $v_{\mf p}(\frac{1}{j_2}) > 0$. 
  The same is true if we interchange $j_1$ and $j_2$. We conclude
  that the set of valuations $v$ of $K_B'$ with $v(j_i) > 0$ does not
  depend on $i$. Therefore, $X_1 = X_2$. 
\end{proof}

\subsection{The scheme of cusps}
\noindent
To describe the boundary of $\ov M^2(f)$, we introduce the {\it scheme
of cusps}, which we call {\it Cusps}. Let $\mf r$ be the intersection 
of all height $1$ primes $\mf p$ containing $\frac{1}{j}$, i.e.,  
$\mf r = {\rm rad}(\frac{1}{j})$. And $V(\mf r) = \ov M^2(f) 
\backslash M^2(f)$. Let $\wh C := \underset{\ll}{\lim}~C/\mf r^n$.

\begin{lemma} \label{lem_Chat}
  The ring $\hat C$ is normal and a finite $A_f[\![\frac{1}{j}]\!]$-algebra. 
\end{lemma}
\begin{proof}
  The ring $B$ is regular. So $C = \oplus_i C_i$ where each 
  $C_i$ is an integrally closed domain. The ring $C$ is excellent. 
  By \cite[7.8.3.vii]{EGA42} it follows that $\wh C$ is normal.
  As $C'$ is a finite $A_f[\frac{1}{j}]$-algebra, it follows that 
  $\wh C$ is a finite $A_f[\![\frac{1}{j}]\!]$-algebra.   
\end{proof}

\noindent
We denote $A_f(\!(\frac{1}{j})\!) := A_f[\![\frac{1}{j}]\!][j]$.
Furthermore, we define the formal scheme
$$\wh {\it Cusps} := {\rm Spf}(\wh C),$$ which is the formal neighbourhood
of {\it Cusps}. Let $\O$ denote the structure 
sheaf of $\wh {\it Cusps}$. The scheme of cusps is defined as
$${\it Cusps} := (\wh {\it Cusps}, \O/ \mf r) = \Spec(C/\mf r).$$

\begin{theorem} \label{thm_extWP}
  The $A_f$-morphism $w_f$ given by Theorem \ref{thm_map}
  can be extended to an $A_f$-morphism $$w_f: \ov M^2(f) \lr M^1(f).$$
  Its restriction to the scheme of cusps gives a finite $A_f$-morphism
  $$w_f: {\it Cusps} \lr M^1(f).$$
\end{theorem}
\begin{proof}
  The Weil pairing gives an $R$-algebra structure $R \lr B$. Because $R$ is
  integral over $A_f$, $B$ is the integral closure of $A_f[j]$ in the quotient
  ring of $B$ and $C$ is the integral closure of $A_f[\frac{1}{j}]$ in this
  quotient ring, it follows immediately, that the ring
  homomorphism $R \lr B$ gives a ring homomorphism $R \lr C$. These
  two maps glue to $$w_f: \ov M^2(f) \lr M^1(f).$$
  The restriction of $w_f$ to ${\it Cusps}$ is given by $R \lr C \lr C/\mf r$. 
  As $C$ is finite over $A_f[\frac{1}{j}]$, it follows that $C/\mf r$ is 
  finite over $A_f$. As $R$ is finite over $A_f$, we may conclude that 
  $w_f$ restricted to ${\it Cusps}$ is finite.
\end{proof}


\section{The cusps and the Tate-Drinfeld module} 
\label{sec_Cusp_TD}

\noindent
In the previous section, we defined the scheme of cusps
and the formal scheme $\wh {\it Cusps} = {\rm Spf}(\wh C)$. In this section,
we will relate these schemes to the universal Tate-Drinfeld module, which we
introduced in Section \ref{sec_UP}. In fact, using the universal property
of the universal Tate-Drinfeld module and the ${\rm Cl}(A) \times
\Gl_2(A/fA)$-equivariance of the Weil-pairing, we will be able to
prove the following theorem. If we write `$\oplus_{\mf p}$', we mean 
the direct sum over all minimal primes $\mf p$ containing $\frac{1}{j}$. 

\begin{theorem} \label{thm_CompToTD}
  There exists an $R[\![\frac{1}{j}]\!]$-linear isomorphism
  $$\wh C \cong \oplus_{\mf p} \lim_{\ll} C/\mf p^n 
  \stackrel{\sim}{\lr} {\oplus} R[\![x]\!]_{(\mf m, \s_i)},$$
  such that $${\it Cusps} \stackrel{\sim}{\lr} \oplus R_{(\mf m, \s_i)}.$$
\end{theorem}

\noindent
From this theorem we can derive the following important corollary:

\begin{cor} \label{cor_cusp}
  The compactification $\ov M^2(f)$ of $M^2(f)$ is regular, and even smooth
  over $\Spec(A_f)$. Furthermore, the scheme of cusps is isomorphic to
  $${\it Cusps} \cong \coprod_{(\mf m, \s_i)} M^1(f),$$ where $\mf m$ runs
  through ${\rm Cl}(A)$ and $\s_i$ runs through the cosets of $N\bs 
  \Gl_2(A/fA)$ where $$N = \left( \begin{array}{cc}
  \F_q^* & A/fA \\ 0 & (A/fA)^* \end{array} \right) \subset \Gl_2(A/fA).$$
  Consequently, the scheme ${\it Cusps}$ consists of $\frac{h(A) \cdot \#
  \Sl_2(A/fA)}{\# (A/fA) \cdot (q-1)}$ copies of $M^1(f)$.
\end{cor}
\begin{proof}
  By Theorem \ref{thm_CompToTD} the ring $C$ is regular in the points
  above $\frac{1}{j}$ and thus $C$ is regular. Consequently, $\ov M^2(f)$
  is regular. The description of ${\it Cusps}$ and the number of
  its components follows from Theorem \ref{thm_CompToTD}.\\
  To prove smoothness over $\Spec(A_f)$, note that by the corollary to
  Proposition 5.4 in \cite{Drin74}, the morphism $M^2(f) \lr \Spec(A_f)$
  is smooth. So we only need to prove smoothness in the closed points
  of ${\it Cusps}$. We have $$\wh C \cong \oplus R[\![x]\!]_{(\mf m, \s_i)},$$
  so $\wh C$ is formally smooth over $A_f$. This implies by 17.5.1 and
  17.5.3 in \cite{EGA44} that the morphism $\ov M^2(f) \lr
  \Spec(A_f)$ is smooth in the closed points of ${\it Cusps}$.
\end{proof}

\subsection{The proof of Theorem \ref{thm_CompToTD}}

\noindent
The rest of this section is devoted to proving Theorem \ref{thm_CompToTD}. 
The universal Tate-Drinfeld module over $\Cal Z$ gives rise
to an $A_f$-morphism $$\Cal Z_{\rm open} \lr M^2(f),$$ where $\Cal Z_{\rm
open}$ denotes the localization of $\Cal Z$ at $(x)$, i.e.,
$$\Cal Z_{\rm open} = \Spec(\oplus R(\!(x)\!)_{(\mf m, \s_i)}).$$
It follows from Remark \ref{rem_Act} that this morphism is ${\rm Cl}(A)\times
\Gl_2(A/fA)$-equivariant.

\par\bigskip\noindent
Let $\Cal Z_{x=0}$ denote the scheme $$\Spec\left(
\oplus R[\![x]\!]_{(\mf m, \s_i)} /(x) \right).$$
The line of argument is as follows: in Lemma \ref{lem_ClAequiv} we
show how to relate the Tate-Drinfeld module to the study of the
cusps, and in Lemma \ref{lem_isoCusps} we describe the 
scheme ${\it Cusps}$. This latter lemma enables us to lift 
the isomorphism ${\it Cusps} \lr \Cal Z_{x = 0}$ to
an isomorphism $\wh C \lr \Cal Z$. Let $\wh C_{\mf p}$ denote the 
completion of the local ring of $C$ at $\mf p$. 

\begin{lemma} \label{lem_ClAequiv}
  Every $R[\![x]\!]_{(\mf m, \s_i)}$ is a finite
  $A_f[\![\frac{1}{j}]\!]$-algebra. Consequently, the morphism $\Cal
  Z_{\rm open} \lr M^2(f)$ comes from a morphism 
  $$h_1: \Cal Z \lr \Spec(C) \subset \ov M^2(f).$$
  Moreover, the universal property of the Tate-Drinfeld module 
  gives rise to a morphism $$h_2: \Spec(\oplus_{\mf p} \wh C_{\mf p}) 
  \lr \Cal Z.$$ The composition $h_1 \circ h_2$ is the natural 
  morphism. 
\end{lemma}
\begin{proof}
  Write $R[\![x]\!] = R[\![x]\!]_{(\mf m, \s_i)}$ for some pair 
  $(\mf m, \s_i)$ and write $(\tdm)_a = \sum_i c_i \t^i$ for the element 
  $a \in A$ which is used to define $j_a$. Then $c_{2\deg(a)}
  \in R(\!(x)\!)^*$, because $\tdm$ is a Drinfeld module over $R(\!(x)\!)$. 
  Moreover, $$\tdm \mod (x) = \psi.$$ So the coefficient $c_{\deg(a)} 
  \in R[\![x]\!]^*$. 
  This implies by definition, that $\frac{1}{j}$ is mapped to
  $\alpha \cdot x^k \in R[\![x]\!]$, with $\alpha \in R[\![x]\!]^*$ and $k \in
  \Z_{>0}$. From this it follows that $R[\![x]\!]$ is a finite
  $R[\![\frac{1}{j}]\!]$-module.
  The morphism $$\Cal Z_{\rm open} \lr M^2(f)$$ comes from a ring homomorphism
  $$C[j] \lr \left( \oplus R[\![x]\!]_{(\mf m, \s_i)} \right) 
  \otimes_{A_f[\![\frac{1}{j}]\!]} A_f(\!(\frac{1}{j})\!),$$ 
  Because $C$ is finite over $A_f[\frac{1}{j}]$ and $R[\![x]\!]$ is 
  finite over $A_f[\![\frac{1}{j}]\!]$, it follows that the image of
  $C$ under this ring homomorphism lies in $\oplus R[\![x]\!]_{(\mf m, 
  \s_i)}$.\\
  For the `moreover'-part, let $\mf p \subset C$ be a minimal prime 
  ideal containing $\frac{1}{j}$. The ring $\wh C_{\mf p}$ 
  is a complete discrete valuation ring and comes equipped with a 
  Tate-Drinfeld structure via the morphism $$\Spec(K_{\wh C_{\mf p}}) 
  \lr M^2(f).$$ By Theorem \ref{thm_univTD}
  there exists a unique ring homomorphism $$\oplus 
  R[\![x]\!]_{(\mf m, \s_i)} \lr \wh C_{\mf p}$$ which induces on 
  $\wh C_{\mf p}$ this Tate-Drinfeld structure.
  This can be done for every minimal prime $\mf p$ containing $\frac{1}{j}$. 
\end{proof}

\begin{lemma} \label{lem_isoCusps}
  The morphism $h_1$ induces an isomorphism 
  $$\Cal Z_{x = 0} \stackrel{\sim}{\lr} {\it Cusps}.$$ 
  Every pair $(\mf m, \s_i)$ corresponds via this isomorphism to one and
  only one minimal prime $\mf p \subset C$ containing $\frac{1}{j}$.
  Consequently, $$\wh C = \oplus_{\mf p} \lim_{\ll} C/\mf p^n.$$ 
\end{lemma}
\begin{proof}
  We will first prove that the number of irreducible components of 
  ${\it Cusps}$ equals the number of irreducible components of 
  $\Cal Z_{x = 0}$. Subsequently, we will show that these components
  of ${\it Cusps}$ intersect nowhere.\\ 
  Because ${\it Cusps} \cong \Spec(C/(\cap \mf p))$,
  the irreducible components of ${\it Cusps}$ are in a one-to-one 
  corresponcende to the minimal primes containing $\frac{1}{j}$.
  The morphisms $h_1$ and $h_2$ introduced in Lemma \ref{lem_ClAequiv} 
  give rise to the following maps on the sets of irreducible components:
  {\small
  $$\{ \mbox{irr. comp. of ${\it Cusps}$} \} \stackrel{h_2}{\lr}
  \{ \mbox{irr. comp. of $\Cal Z_{x = 0}$} \} \stackrel{h_1}{\lr}
  \{ \mbox{irr. comp. of ${\it Cusps}$} \}.$$}

  \noindent
  As a morphism of schemes $h_1 \circ h_2$ is the natural map. Therefore, 
  the composition $h_1 \circ h_2$  
  on the set of irreducible components is the identity. Consequently, 
  $h_2$ on the irreducible components is injective.\\
  Moreover, the set of irreducible components of $\Cal Z_{x = 0}$ is by 
  definition one orbit under the elements $(\mf m, \s_i)$. 
  Clearly, the map $h_1$ on the set of connected components is equivariant 
  under this group action, 
  and as the image of $h_2$ is not empty, it follows that the first
  map on the irreducible components is also surjective. So we may conclude
  that the number of irreducible components of ${\it Cusps}$ equals the
  number of irreducible components of $\Cal Z_{x = 0}$.\\
  The irreducible components of $\Cal Z_{x = 0}$ intersect nowhere.
  We will prove that this is also the case for the irreducible components of
  ${\it Cusps}$. By the extension of the Weil pairing to
  $\ov M^2(f)$, the ring $C$ comes equipped with an $R$-algebra structure.
  Let $$\zeta: R \lr \oplus R[\![x]\!]_{(\mf m, \s_i)}$$ denote the 
  composition $R \stackrel{w_f^{\#}}{\lr} C \stackrel{h_1^{\#}}{\lr} 
  \oplus R[\![x]\!]_{(\mf m, \s_i)}.$\\
  Choose any maximal ideal $\mf n \subset R$ and let $\mf q$ run over all
  minimal primes of $C/\mf n C$ containing $\frac{1}{j}$. We write
  $$\wh{C/\mf n C} = \lim_{\ll} (C/\mf n C)/ (\frac{1}{j})^n,$$
  and $\wh{(C/\mf n C)}_{\mf q}$ for the completion along $\mf q$ of
  the local ring $(C/\mf n C)_{\mf q}$.\\
  In this case we have analogues of the morphisms $h_1, h_2$, namely, 
  $R$-algebra homomorphisms
  $$\begin{CD} \wh{C/\mf n C} @> \tilde h_1 >>
  \oplus R/\zeta(\mf n)[\![x]\!]_{(\mf m, \s_i)} @> \tilde h_2 >>
  \oplus \wh{C/\mf n C}_{\mf q}.
  \end{CD}$$
  And now, as before, we define maps on the sets of irreducible components:
  {\small
  $$\tilde h_2: \{ \mbox{irr. comp. of ${\it Cusps} \times \Spec(R/\mf n)$}
  \}
  \lr \{ \mbox{irr. comp. of $\Cal Z_{x = 0} \times \Spec(R/\mf n)$} \},$$ 
  $$\tilde h_1: \{ \mbox{irr. comp. of $\Cal Z_{x = 0}\times
  \Spec(R/\mf n)$} \} \lr \{ \mbox{irr. comp. of
  ${\it Cusps}\times \Spec(R/\mf n)$} \}.$$}

  \noindent
  The composition of these two maps on the set of irreducible components
  is the identity. Namely, the composition $\tilde h_1 \circ
  \tilde h_2$ is the natural map on the rings. Using the same argument
  as above shows that $\tilde h_2$ on the irreducible components is a 
  bijection.\\
  We conclude that for every prime $\mf n \subset R$ the number of
  irreducible components of $${\it Cusps} \times \Spec(R/\mf n)$$ equals the
  number of irreducible components of ${\it Cusps}$. Recall that by Theorem 
  \ref{thm_extWP} the morphism $w_f: {\it Cusps} \lr M^1(f)$ is finite. 
  Therefore, if the irreducible components
  would intersect above some prime ideal $\mf n \subset R$, then
  ${\it Cusps}\times \Spec(R/\mf n)$ would have less irreducible components.
  As this is not the case, we conclude that the irreducible components
  of ${\it Cusps}$ intersect nowhere.\\
  We write $${\it Cusps} = \Spec(\oplus S_{(\mf m, \s_i)})$$ where
  $\Spec(S_{(\mf m, \s_i)})$ are the connected components of ${\it Cusps}$. 
  For every
  pair $(\mf m, \s_i)$, we get $R$-linear ring homomorphisms on the
  connected components:
  $$S_{(\mf m, \s_i)} \stackrel{h_1^{\#}}{\lr} R \stackrel{h_2^{\#}}
  \lr S_{(\mf m, \s_i)}.$$ 
  Because the composition is the identity and $S_{(\mf m, \s_i)}$ is 
  finite over $R$, it follows that $$S_{(\mf m, \s_i)} \cong R.$$
  For the latter two statements of the lemma, note that the isomorphism
  implies that the minimal primes $\mf p$ are relatively prime and,
  consequently, $\mf r = \prod \mf p.$
\end{proof}

\noindent
The next step is to lift the isomorphism from Lemma 
\ref{lem_isoCusps} to an isomorphism $$\wh C \stackrel{\sim}{\lr}
\oplus R[\![x]\!]_{(\mf m, \s_i)}.$$ Let
$$\Cal W := \underset{\ll}{\lim}~ C/\mf p^n$$ for some minimal prime
$\mf p$ containing $\frac{1}{j}$. 

\begin{lemma} \label{lem_isoP}
  The ring $\Cal W$ is isomorphic to $R[\![x]\!]$.
\end{lemma}
\begin{proof}
  By Lemma \ref{lem_Chat} the ring $\Cal W$ is integrally closed and a finite 
  $A_f[\![\frac{1}{j}]\!]$-algebra, and by the isomorphism of Lemma
  \ref{lem_isoCusps}, the ring $\Cal W$ is a finite 
  $R[\![\frac{1}{j}]\!]$-algebra.\\
  The completion of the local ring $\Cal W_{\mf p}$ is isomorphic
  to $\wh C_{\mf p}$. The morphisms $h_1$ and $h_2$ give on the 
  completions of the local rings injective maps 
  $$\wh C_{\mf p} \stackrel{h_1^{\#}}{\lr} K_R[\![x]\!]
  \stackrel{h_2^{\#}}{\lr} \wh C_{\mf p}.$$
  As $h_2^{\#} \circ h_1^{\#}$ is the identity, there exists an isomorphism 
  $K_R[\![x]\!] \cong \wh C_{\mf p}.$\\ 
  We conclude that $\Cal W$ is regular. Therefore, we may asssume that   
  $x$ is an element of $\Cal W$ and that $\Cal W$ is a finite 
  $R[\![x]\!]$-algebra. 
  So we get injective $R[\![x]\!]$-linear ring homomorphisms 
  $$R[\![x]\!] \lr \Cal W \lr R[\![x]\!]$$ where the first
  map is the $R[\![x]\!]$-structure morphism of $\Cal W$ and the second 
  map is $h_1^{\#}$. We conclude that $\Cal W \cong R[\![x]\!]$.
\end{proof}

\noindent
This enables us to prove Theorem \ref{thm_CompToTD}:

\begin{proof}[Proof of Theorem \ref{thm_CompToTD}] 
  By the previous lemma, it follows that $$\wh C \cong \oplus
  R[\![x]\!]_{(\mf m, \s_i)}.$$ Together with Lemma \ref{lem_isoCusps}
  the theorem follows.
\end{proof}


\section{Components of $M^2(f)$}
\label{sec_compon}

\noindent
In this section we describe the geometric components of $\ov M^2(f)$ and prove
the connectedness of $M^2(f)$. For a non-zero prime $\mf P \subset R$ we
write $\kappa(\mf P) := R/\mf P$. The first result is the following:

\begin{theorem} \label{thm_geom_comp}
  The scheme $$\ov M^2(f) \underset{A_f}{\times} M^1(f)$$
  consists of $h(A)\cdot [(A/fA)^*:\F_q^*]$ connected components, which
  are all geometrically connected. Moreover, for every non-zero prime 
  ideal $\mf P \subset R$ the fibre at $\mf P$
  $$\ov M^2(f) \underset{A_f}{\times} \Spec(\kappa(\mf P)),$$
  consists of $h(A) \cdot [(A/fA)^*: \F_q^*]$ connected components, which are
  all geometrically connected. 
\end{theorem}
\begin{proof}
  Let $K_{\infty}$ be the completion of the quotient field of $A_f$
  along the point $\infty$, and let $\C_{\infty}$ denote the
  completion of the algebraic closure of $K_{\infty}$. By the analytic
  theory, as is shown in \cite{PT1}, we know that $$\ov M^2(f)
  \times_{A_f} \Spec(\C_{\infty})$$ consists of $h(A) \cdot
  [(A/fA)^*:\F_q^*]$ components. Because $R$ is a Galois extension of
  $A_f$ with Galois group $G$, we have $R\otimes_{A_f} R \cong
  \oplus_G R.$ By the Weil pairing one sees that $$\ov M^2(f)
  \underset{A_f}{\times} M^1(f) \stackrel{w_f}{\lr} M^1(f)
  \times_{A_f} M^1(f)$$ consists of $h(A)\cdot \# G$ connected components.
  As $\# G = [(A/fA)^*:\F_q^*]$, these components are geometrically
  connected components.\\
  Consider the fibres over $R$. Let $\mf P \subset R$ be a non-zero prime 
  ideal and let $V$ be the
  completion along $\mf P$ of the local ring $R_{\mf P}$. Suppose
  $\ov M^2(f) \times_R \Spec(\kappa(\mf P))$ has more than one
  connected component, then also $\ov M^2(f) \times_R \Spec(V/\mf
  P^n)$ has more than one connected component for every $n$ and
  consequently, both $\ov M^2(f) \times \Spec(V)$ and $\ov M^2(f)
  \times_R \Spec(K_V)$ consist of more than one component. This,
  however, contradicts the fact that $\ov M^2(f) \times_R M^1(f)$ is
  geometrically connected. So we conclude that
  $\ov M^2(f) \times_R \Spec(\kappa(\mf P))$ is geometrically connected.
\end{proof}

\noindent
This theorem enables us to say something about the Drinfeld modular curves.
Let, as before, $N = \left( \begin{array}{cc} \F_q^* & A/fA \\ 0 & (A/fA)^*
\end{array} \right)$.

\begin{theorem}
  For every $R$-field $K$ the curve $\ov M^2(f) \times_R \Spec(K)$ is
  a smooth, irreducible curve containing $h(A) \cdot [\Gl_2(A/fA): N]$
  cusps.
\end{theorem}
\begin{proof}
  Clearly, the scheme {\it Cusps} consists of $h(A) \cdot [\Gl_2(A/fA): N]$
  copies of $R$. Consequently, $${\it Cusps} \times \Spec(K)$$ consists
  of $h(A) \cdot [\Gl_2(A/fA): N]$ points. The irreducibility follows
  immediately from the proof of Theorem \ref{thm_geom_comp}.
\end{proof}

\subsection{The analogue of $X_0(N)$}
\noindent
The analogue in the setting of Drinfeld modular curves of the modular
curve $X_0(N)$ is the curve $$X_0(f):= \ov M^2(f)/H, \quad {\rm where}~ 
H = \left( \begin{array}{cc}
(A/fA)^* & A/fA\\ 0 & (A/fA)^* \end{array} \right) \subset \Gl_2(A/fA).$$
One may deduce from Theorem \ref{thm_geom_comp} the following
theorem concerning the cusps and the geometric components of $\ov M^2(f)/H$.
Define $R_0 = R^{(A/fA)^*/\F_q^*}$, i.e., $\Spec(R_0) = M^1(1)$. 
Write ${\it Cusps}_0$ for the scheme of cusps of $X_0(f)$.

\begin{theorem}
  The Weil pairing induces an isomorphism 
  $${\it Cusps}_0 \stackrel{\sim}{\lr} \coprod_{(\mf m, \rho)} M^1(1)$$ where 
  $\rho$ runs through the double cosets $N \bs \Gl_2(A/fA) / H$.\\
  The scheme $X_0(f)$ is connected, and for any $R_0$-field $K$ the 
  scheme $X_0(f) \times \Spec(K)$ consists of $h(A)$ geometrically
  connected components. 
\end{theorem}
\begin{proof}
  The morphism $w_f$ gives an isomorphism between $M^1(f)$ and
  any connected component of ${\it Cusps}$. Consequently, $w_f$ gives
  an isomorphism ${\it Cusps} \lr \oplus_{(\mf m, \s_i)} M^1(f)$. 
  Recall that $N$ acts trivially on each copy of $M^1(f)$. Furthermore,
  the action of $\s = \left( \begin{array}{cc} \alpha & 0 \\ 0 & 
  \alpha \end{array} \right)$ with $\alpha \in (A/fA)^*$ on ${\it Cusps}$
  is as follows. Let $i,j\in \N$ such that $\s_i \circ \s \in N \s_j$, 
  and consider $\alpha$ as an element of the Galois group 
  $\Gal(K_R/K_{R_0}) \cong (A/fA)^*/\F_q^*$,
  then $\s$ acts as $$M^1(f)_{(\mf m, \s_i)} \stackrel{\alpha}{\lr} 
  M^1(f)_{(\mf m, \s_j)}.$$ 
  Because $H$ contains the subgroup $\left\{ \left( \begin{array}{cc} \alpha 
  & 0 \\ 0 & \alpha \end{array} \right) \mid \alpha \in (A/fA)^* \right\}$, 
  we see that dividing out ${\it Cusps}$ by $H$ gives an isomorphism
  $${\it Cusps}_0 \stackrel{\sim}{\lr}  
  \oplus_{(\mf m, \rho)} M^1(1)$$ where 
  $\rho$ runs through the cosets of $N\bs \Gl_2(A/fA) / H$.\\ 
  The number of components follows immediately from Theorem 
  \ref{thm_geom_comp}. 
\end{proof}


\def\cprime{$'$} \def\cprime{$'$}
  \def\polhk#1{\setbox0=\hbox{#1}{\ooalign{\hidewidth
  \lower1.5ex\hbox{`}\hidewidth\crcr\unhbox0}}}

\end{document}